\documentclass[12pt,a4paper,svgnames]{amsart}

\usepackage{amsthm,amsmath,amsfonts,amssymb}
\usepackage{hyperref}
\usepackage[usenames]{color}
\usepackage{xcolor}
\usepackage{mathtools}
\usepackage{amsmath,amsfonts}
\usepackage{calrsfs}
\usepackage{graphicx}
\usepackage{tikz}
\usepackage{dsfont}
 \usepackage{mathrsfs}
 \usetikzlibrary{arrows}

\renewcommand{\proof}{{\noindent \bf Proof: }}
\newtheorem{theorem}{Theorem}[section]
\newtheorem{lemma}[theorem]{Lemma}
\newtheorem{corollary}[theorem]{Corollary}
\newtheorem{definition}[theorem]{Definition}

\newtheorem{remark}{Remark}[section]

\title[Fractional oscillon equation]{Fractional oscillon equations; solvability and connection with classical oscillon equations}

\numberwithin{equation}{section} \numberwithin{theorem}{section}

\author[F. D. M. Bezerra]{Flank D. M. Bezerra}
\address[F. D. M. Bezerra]{Universidade Federal da Para\'{\i}ba, Departamento de Matem\'atica, 58051-900 Jo\~{a}o Pessoa PB, Brazil.}
\email{flank@mat.ufpb.br}

\author[R. N. Figueroa-L\'opez]{Rodiak N. Figueroa-L\'opez$^{\star\star}$}\thanks{$^{\star\star}$ Research partially supported by CAPES/Finance Code 001/2019, Brazil.}
\address[R. N. Figueroa-L\'opez]{Universidade Federal de S\~{a}o
Carlos, Departamento de Matem\'atica, 13565-905 S\~{a}o
Carlos SP, Brazil.}
\email{rodiak@dm.ufscar.br}

\author[M. J. D. Nascimento]{Marcelo J. D. Nascimento$^\star$}\thanks{$^\star$Research partially
supported by FAPESP \#2017/06582-2, Brazil}
\address[M. J. D. Nascimento]{Universidade Federal de S\~{a}o
Carlos, Departamento de Matem\'atica, 13565-905 S\~{a}o
Carlos SP, Brazil.}
\email{marcelo@dm.ufscar.br}
\date{\today}

\begin{document}

\begin{abstract} 
In this paper we are concerned with the asymptotic behavior of nonautonomous fractional approximations of oscillon equation
\[
u_{tt}-\mu(t)\Delta u+\omega(t)u_t=f(u),\ x\in\Omega,\ t\in\mathbb{R},
\]
subject to Dirichlet boundary condition on $\partial \Omega$, where $\Omega$ is a bounded smooth domain in $\mathbb{R}^N$, $N\geqslant 3$, the function $\omega$ is a time-dependent damping, $\mu$ is a time-dependent squared speed of propagation, and $f$ is a nonlinear functional. Under structural assumptions on $\omega$ and $\mu$ we establish the existence of  time-dependent attractor for the fractional models in the sense of Carvalho, Langa, Robinson \cite{CLR}, and Di Plinio, Duane, Temam \cite{DDT1}.

\vskip .1 in \noindent {\it Mathematical Subject Classification 2010:} 37B55, 35B40, 35B41, 34A08, 35L71.
\newline {\it Key words and phrases:} oscillon equation; fractional powers; fractional equations; pullback attractor.
\end{abstract}

% 	37B55   	Nonautonomous dynamical systems
% 35B40   	Asymptotic behavior of solutions
% 35B41   	Attractors
 %	34A08   	Fractional differential equations
  %35L71   	Semilinear second-order hyperbolic equations

\maketitle
\tableofcontents

\section{Introduction}

In this paper we  consider the oscillon-type equations  that can be applied in cosmology to model and represent some transient persistent structures. Theses equations has been analyzed in literature in the sense of attractors for nonlinear evolution processes; namely,  in Di Plinio, Duane, Temam \cite{DDT1} the authors consider the equation
\begin{equation}\label{asdar}
u_{tt}-e^{-2Ht}\Delta u+Hu_t=f(u),\ x\in(0,1),\ t>\tau,\ \tau\in\mathbb{R},\\
\end{equation}
and the asymptotic behavior of solutions in considered with periodic boundary conditions, $H>0$ and $f$ is defined by a potential of arbitrary polynomial growth. The authors shows the existence of  regular global attractor, and proves that the  kernel section have finite fractal dimension.

In Di Plinio, Duane, Temam \cite{DDT2} the authors consider the equation
\begin{equation}\label{asdar01}
u_{tt}-\mu(t)\Delta u+\omega(t)u_t=f(u),\ x\in\Omega,\ t>\tau,\ \tau\in\mathbb{R},\\
\end{equation}
and the asymptotic behavior of solutions in considered with zero Dirichlet boundary condition on $\partial\Omega$, where $\Omega$ is a bounded smooth domain in $\mathbb{R}^3$. In this paper, the authors shows the existence of pullback attractor and a result of regularity of the attractors as well the finiteness of fractal dimension of the kernel section.

In Conti, Pata and Temam \cite{CPT} results of the recent theory of attractors in time-dependent spaces from \cite{DDT1} are generalized.

The aim of this paper is to analyze nonlinear evolution processes defined on time-dependent spaces generated by a class of $N$-dimensional fractional oscillon equations. More precisely, we consider the following initial-boundary value problem
\begin{equation}\label{xxx}
\begin{cases}
u_{tt}+\omega(t)u_t-\mu(t)\Delta u=f(u),& x\in\Omega,\ t>\tau,\\
u(x,t)=0,& x\in\partial\Omega,\ t\geqslant\tau,\\
u(x,\tau)=u_\tau(x), u_t(x,\tau)=v_\tau(x),& x\in\Omega,
\end{cases}
\end{equation}
where $\Omega$ is a bounded smooth domain in $\mathbb{R}^N$, $N\geqslant  3$, and  $f \in C^{1}(\mathbb{R})$ satisfies
\begin{equation}\label{growth-condition}
|f'(s)| \leqslant  C(1+ |s|^{\rho -1}),\quad s\in \mathbb{R},
\end{equation}
for some
\begin{equation}\label{eq:supercritical}
\frac{N}{N-2}< \rho < \frac{N+2}{N-2},
\end{equation}
and
\begin{equation}\label{CondDissip}
\limsup_{|s|\to \infty} \frac{f(s)}{s} \leqslant 0.
\end{equation}

The damping coefficient $\omega:\mathbb{R}\to\mathbb{R}^+$ is assumed to be decreasing strictly positive diffe\-rentiable function, with $\omega(t)$ bounded as $t\to -\infty$ and we set 
\begin{equation}\label{132r}
W:=\sup_{t\in\mathbb{R}}\omega(t)<\infty.
\end{equation}
The degeneracy $\displaystyle\lim_{t\to\infty}\omega(t)=0$ is allowed.  Suppose also that the function $\omega$ is $(\zeta,\kappa_0)-$H\"{o}lder continuous in $\mathbb{R}$; that is, 
\begin{equation}\label{hol-omega}
 |\omega(t)-\omega(\tau)|\leqslant\kappa_0|t-\tau|^{\zeta},\quad \forall t,\tau\in\mathbb{R},
\end{equation}
where $\zeta\in(0,1]$ and $\kappa_0>0$.
Concerning to $\omega$ we define the {\it decay rate} $\varepsilon_{\omega}:\mathbb{R}\to\mathbb{R}^+$, given by 
\begin{equation}\label{condE}
\varepsilon_{\omega}(t)=\min\Big\{1,\frac{\omega(t)}{4},\frac{c_1}{4(W+2)},\frac{d^{2}_{0}}{3d^{2}_{1}}\Big\},
\end{equation}
where $c_1,d_0,d_1> 0$ are positive constants will to be specified later.

The main structural assumption on $\mu$ is as follows: there are positive constants $\mu_{min}$ and $\mu_{max}$ such that
\begin{equation}\label{123}
0<\mu_{min}\leqslant\mu(t)\leqslant\mu_{max}, \quad\forall t\in\mathbb{R}.
\end{equation}
Additionally, suppose that there exists a function $\vartheta:\mathbb{R}\to(0,\infty)$ satisfying
\begin{equation}\label{1506r}
0<\mu^{\prime}(t)\leqslant\vartheta(t)\mu(t)\quad\text{ with } \quad\sup_{\tau\leqslant t}\vartheta(t)\leqslant \varepsilon^{2}_{\omega}(t),\quad\forall t\in\mathbb{R}.
\end{equation}
Furthermore, suppose that the function $\mu$ is $(\gamma,\kappa)-$H\"{o}lder continuous in $\mathbb{R}$; that is, 
\[
|\mu(t)-\mu(\tau)|\leqslant\kappa|t-\tau|^{\gamma}, \quad\forall t,\tau\in\mathbb{R},
\]
where $\gamma\in[1/2,1]$ and $\kappa>0$.

 We first introduce some notations and terminologies needed to better understanding our results. Let $X= L^2(\Omega)$ equipped with the standard norm $\|\cdot\|$ and scalar product $\langle\cdot,\cdot\rangle$, and let $A:D(A)\subset X\to X$ be the unbounded linear operator 
\begin{equation}\label{eq:A}
Au=-\Delta u\ \text{ for } \ u\in D(A) = H^2(\Omega) \cap H^{1}_{0}(\Omega).
\end{equation}
 Is is well known that $A$ is a positive, self-adjoint operator and $-A$ generates a compact analytic $C^0$-semigroup on $X$. Denote by $X^\alpha$ the fractional power spaces associated with the operator $A$; that is, $X^\alpha = D(A^\alpha)$ with the norm $\|A^\alpha \cdot \|:X^\alpha\to \mathbb{R}^+$, and this norm is induced by the inner product
\[
\langle u,v \rangle_{X^\alpha}=\langle A^\alpha u, A^\alpha v \rangle.
\] 
For $\alpha>0$, define $X^{-\alpha}$ as the completion of $X$ with the norm $\|A^{-\alpha}\cdot\|$, (see Amann \cite{Am}). With this notation we have $X^\frac{1}{2}=H^1_0(\Omega)$ and $X^1=H^2(\Omega)\cap H^1_0(\Omega)$. It follows from Amann \cite[Chapter V]{Am} that $X^{-\alpha}=(X^{\alpha})'$ for any $\alpha>0$.

Consider the Hilbert space given by
\[
Y_t=X^{\frac{1}{2}}\times X
\]
equipped with the inner product
\[
\left\langle \left[\begin{smallmatrix} u_1\\ v_1 \end{smallmatrix}\right], \left[\begin{smallmatrix} u_2\\ v_2 \end{smallmatrix}\right] \right\rangle_{Y_t}=\mu(t)^{\frac{1}{2}}\langle u_1, u_2 \rangle_{X^\frac{1}{2}}+\langle u_1, u_2 \rangle+\langle v_1, v_2 \rangle.
\]

Note that the spaces $Y_t$ are all the same as linear spaces and the norms $\|\cdot\|_{Y_t}^2=\langle\cdot,\cdot\rangle_{Y_t}$ and $\|\cdot\|_{Y_\tau}^2=\langle\cdot,\cdot\rangle_{Y_\tau}$ are equivalents for any fixed $t,\tau\in\mathbb{R}$.

Let $v=u_t$, where $u$ is the unknown function in \eqref{xxx}. Thus, the problem \eqref{xxx} can be rewrite as a Cauchy problem in $Y_t$; namely
\begin{equation}\label{yyy}
\begin{cases}
\dfrac{dw}{dt}+\varLambda(t)w=F(t,w),\quad t>\tau,\\
w(\tau)=w_\tau,
\end{cases}
\end{equation}
where $w=\left[\begin{smallmatrix} u\\ v \end{smallmatrix}\right]$,  $w_\tau=\left[\begin{smallmatrix} u_\tau\\ v_\tau \end{smallmatrix}\right]$, $\varLambda(t):D(\varLambda(t))\subset Y_t\to Y_t$ is the unbounded linear operator defined by
\[
D(\varLambda(t))=X^{1}\times X^{\frac{1}{2}}=Y_t^1
\]
and 
\begin{equation}\label{eq:Lambda}
\varLambda(t)\left[\begin{smallmatrix} u\\ v \end{smallmatrix}\right]=\left[\begin{smallmatrix} 0 & -I\\ \mu(t)A & 0 \end{smallmatrix}\right]\left[\begin{smallmatrix} u\\ v \end{smallmatrix}\right]:=\left[\begin{smallmatrix} -v\\ \mu(t)Au \end{smallmatrix}\right],\quad \forall \left[\begin{smallmatrix} u\\ v \end{smallmatrix}\right]\in Y_t^1,
\end{equation}
and the nonlinearity $F$ is defined by 
\begin{equation}\label{eq:Lambdaf}
F\left(t,\left[\begin{smallmatrix} u\\ v \end{smallmatrix}\right]\right)=\left[\begin{smallmatrix} 0\\ f^e(u)-\omega(t)v \end{smallmatrix}\right],
\end{equation}
where $f^e:X^{\frac{1}{2}}\to X$ is the Nemikist\u{\i}i operator associated with $f$; namely, $f^e(u)(x):=f(u(x))$ for any $u\in X^{\frac{1}{2}}$ and $x\in\Omega$. By simplicity of notation, we will denote $f^e$ simply by $f$.

Since \eqref{xxx} can be viewed in the form \eqref{yyy}, due to spectral properties of the operator $\varLambda(t)$, as will see below, it is natural to consider with \eqref{yyy} a family of Cauchy problems in $Y_t$
\begin{equation}\label{ffhghgf}
\begin{cases}
\dfrac{dw^\alpha}{dt}+\varLambda(t)^\alpha w^\alpha=F(t,w^\alpha),\quad t>\tau,\\
w^\alpha(\tau)=w_{\tau\alpha},
\end{cases}
\end{equation}
where $0<\alpha<1$, $w^\alpha=\left[\begin{smallmatrix} u^\alpha\\ v^\alpha \end{smallmatrix}\right]$,  $w_{\tau\alpha}=\left[\begin{smallmatrix} u_{\tau\alpha}\\ v_{\tau\alpha} \end{smallmatrix}\right]$, and the unbounded linear operator $\varLambda(t)^\alpha:D(\varLambda(t)^\alpha)\subset Y_t\to Y_t$ is the fractional power of $\varLambda(t)$  characterized by \eqref{pot-a} bellow. The nonlinearity $F$ is defined by  \eqref{eq:Lambdaf}. 

The main contributions of this paper for the specialized literature can be summarized in three parts:
\begin{itemize}
\item We present approximations of the classical oscillon equation by parabolic type pro\-blems of ``lower''  order  via fractional power;
\item Exploiting parabolic structure of \eqref{ffhghgf} we prove local well posedness of \eqref{ffhghgf} for all $\alpha<1$ close enough to $1$ in a suitably large phase space of initial data containing the energy space $H^1_0(\Omega)\times L^2(\Omega)$;
\item Using \eqref{CondDissip} and exploiting gradient structure of the fractional model \eqref{ffhghgf} we esta\-blish the global well posedness and the existence of pullback attractor, in the sense of Carvalho, Langa, Robinson \cite{CLR}, and existence of time-dependent global attractor, in the sense of Di Plinio, Duane, Temam \cite{DDT1}, for the problem \eqref{ffhghgf}.
\end{itemize}

This paper is organized as follows. In Section \ref{Sec:attractors} we introduce notations, terminologies and results of the theory of attractors in time-dependent spaces.  In Section \ref{Sec:LA} we study the properties of the fractional powers of operators $A$ and $\varLambda(t)$ defined in \eqref{eq:A} and \eqref{eq:Lambda}, respectively. Finally, in Section \ref{Sec:frac-osc-eq} we deduce an $N$-dimentional fractional oscillon equation and we study the gradient structure and dissipativity property of this fractional model.

\noindent {\bf Terminology}. In this paper we will use the term \textit{PDE} to refer to \textit{partial differential equation}.

\section{Attractors in time-dependent spaces}\label{Sec:attractors}

 In this section we recall a few basic definitions and main results concerning with pullback attractors theory in time-dependent spaces, see Conti, Pata and Temam \cite{CPT} and Di Plinio, Duane, Temam \cite{DDT1,DDT2}. 

\begin{definition}
For $t\in\mathbb{R}$, let $X_t$ be a  Banach space endowed with norm $\|\cdot\|_{X_t}$ (see \eqref{3.8} for an example). A (continuous) process $S(\cdot,\cdot)$ is a two-parameter family of mappings $\{S(t,\tau): X_\tau \to X_t: t\geqslant\tau\in\mathbb{R}\}$  with properties:
\begin{itemize}
\item[(i)] $S(t,t)=I_{X_t}$;

\item[(ii)] $S(t,\tau) \in C(X_\tau, X_t )$;

\item[(iii)] $S(t,r)S(r,\tau)=S(t,\tau)$ for $\tau\leqslant r\leqslant t$.
\end{itemize}
\end{definition}

\begin{definition}
A family of subsets $\mathscr{B}=\{\mathscr{B}(t)\subset X_t: t\in\mathbb{R}\}$ is pullback-bounded if
$$
R(t)=\sup_{\tau\in(-\infty,t]}\|\mathscr{B}(\tau)\|_{X_\tau}=\sup_{\tau\in(-\infty,t]}\sup_{z\in \mathscr{B}(\tau)}\|z\|_{X_\tau}<\infty,\quad\forall t\in\mathbb{R}.
$$
\end{definition}

\begin{definition}
A pullback-bounded family $\mathbb{A}=\{\mathbb{A}(t):t\in\mathbb{R}\}$ is called pullback absorbing if for every pullback-bounded family $\mathscr{B}=\{\mathscr{B}(t)\subset X_t: t\in\mathbb{R}\}$ and for every $t\in\mathbb{R}$ there exists $t_0=t_0(t)\leqslant t$ such that
$$
S(t,\tau)\mathscr{B}(\tau)\subset\mathbb{A}(t),\quad \forall \tau\leqslant t_0.
$$
\end{definition}

\begin{definition}
Given a family of sets $\mathscr{B}$, its time-dependent $\omega-$limit is the family $\omega_{\mathscr{B}}=\{\omega_{\mathscr{B}}(t)\subset X_t:t\in\mathbb{R}\}$, where $\omega_{\mathscr{B}}(t)$ is defined as
$$\omega_{\mathscr{B}}(t)=\bigcap_{r\leqslant t}\overline{\bigcup_{\tau\leqslant r }S(t,\tau)\mathscr{B}(\tau)},$$
and the above closures are taken in $X_t$. A more concrete characterization is the following:
$$\omega_{\mathscr{B}}(t)=\{z\in X_t: \exists \tau_n\to -\infty, z_n\in \mathscr{B}(\tau_n) \text{ with } \|S(t,\tau_n)z_n-z\|_{X_t}\to 0\text{ as }n\to\infty\}.$$
\end{definition}

\begin{definition}\label{def:attractor}
A family of compact subsets  $\mathscr{A}=\{\mathscr{A}(t)\subset X_t: t\in\mathbb{R}\}$ is called  time-dependent global attractor for the process $S(\cdot,\cdot)$ if it fulfills the following properties:
\begin{itemize}
\item[$(i)$] (Invariance) $S(t,\tau)\mathscr{A}(\tau)=\mathscr{A}(t)$, for every $\tau\leqslant t$;

\item[$(ii)$] (Pullback attraction) for every pullback-bounded family $\mathscr{B}$ and every $t\in\mathbb{R}$,
$$\lim_{\tau\to-\infty}dist_{X_t}(S(t,\tau)\mathscr{B}(\tau),\mathscr{A}(t))=0,$$
where $dist_{X_t}(\cdot,\cdot)$ is the Hausdorff semidistance in $X_t$.
\end{itemize}
\end{definition}

\begin{remark}\label{remark:2.1}
\noindent
$1)$
To ensure the uniqueness of the pullback attractor in the sense of Definition \ref{def:attractor}, it is necessary to add the condition (see Di Plinio, Duane, Temam \cite{DDT1,DDT2})

$(iii)$ $\mathscr{A}$ is a pullback-bounded family.

\noindent 
$2)$ In the sense of Carvalho, Langa, Robinson \cite[Definition 1.12]{CLR}, condition additional for uniqueness of pullback attractor is $\mathscr{A}$ to be the minimal family of closed sets with property of pullback attraction.

\end{remark}

We will use the Kuratowski measure of noncompactness: if $X$ is a Banach space and $D\subset X$, then
\[
\alpha(D)=\inf\{\delta> 0: D\text{ has a finite cover of balls of }X\text{ of radius less than }\delta\}.
\]

Recall some properties of the Kuratowski measure, redirecting to \cite{hale} for more details and proofs:
\begin{itemize}
\item [(K.1)] $\alpha(D) = 0$ if and only if $D$ is compact in $X$;
\item [(K.2)]$D_1\subset D_2$ implies $\alpha(D_1)\leqslant \alpha(D_2)$; 
\item [(K.3)] $\alpha(D)=\alpha(\overline{D})$, where $\overline{D}$ denotes the closure of the set $D$.
\end{itemize}

\begin{remark}
The shorthand $\alpha_t$ stands for the Kuratowski measure in the space $X_t$. We remark that, for fixed $\tau,t\in\mathbb{R}$, $\alpha_\tau$ and $\alpha_t$ are equivalent measures of noncompactness whenever there is a Banach space isomorphism between $X_\tau$ and $X_t$.
\end{remark}

The proof of the following result can be found in Di Plinio, Duane, Temam \cite[Theorem 2.1]{DDT1}.

\begin{theorem}\label{existence-attractor}
Assume that the nonlinear evolution process $S(\cdot,\cdot)$ possesses a pullback absorbing family $\mathbb{A}$ for which
\begin{equation}\label{measurealpha}
\lim_{\tau\to-\infty}\alpha_t(S(t,\tau)\mathbb{A}(\tau))=0,\quad\forall t\in\mathbb{R},
\end{equation}
where $\alpha_t$ is the Kuratowski measure in the space $X_t$. Then $\mathscr{A}(\cdot)=\omega_{\mathbb{A}}(\cdot)$  is a pullback attractor for $S(\cdot,\cdot)$. 
\end{theorem}

\begin{corollary}\label{uniqueness-attractor}
Under the same assumptions of Theorem \ref{existence-attractor}, we have
\[
\mathscr{A}(t)=\omega_{\mathbb{A}}(t)\subset \mathbb{A}(t),\quad \forall t\in\mathbb{R}.
\]
In particular, $\mathscr{A}$ is a pullback-bounded family, and therefore unique in the sense of $1)$ of Remark \ref{remark:2.1}.
\end{corollary}

\section{Fractional powers of operators}\label{Sec:LA}

In this section we present results on  fractional powers of the operators $A$ and $\varLambda(t)$ defined  in \eqref{eq:A} and \eqref{eq:Lambda}, respectively.

\begin{lemma}\label{lem:useful:properties}
If $A$ and $\varLambda(t)$ are as in \eqref{eq:A} and \eqref{eq:Lambda}, respectively, then we have all the following:

\noindent
i) $0\in\rho(\varLambda(t))$ and
$$
\varLambda(t)^{-1}=\begin{bmatrix} 0 & \mu(t)^{-1}A^{-1}\\ -I & 0
\end{bmatrix}.
$$

\noindent
ii) The adjoint $\varLambda(t)^{*}$ of $\varLambda(t)$ is given by
$$
\varLambda(t)^{*} = \begin{bmatrix} 0 & I\\ -\mu(t)A & 0
\end{bmatrix}=-\varLambda(t).
$$

\noindent
iii) The unbounded linear operator $i\varLambda(t)$ is self-adjoint and
$\varLambda(t)$ is the infinitesimal generator of a $C_0$-group $\{e^{\varLambda(t) \tau} :\tau\geqslant 0\}$ of unitary operators in $Y_t$.

\noindent
iv) Fractional powers $\varLambda(t)^\alpha$ can be defined for $\alpha\in(0,1)$ through
\begin{equation}\label{eq:to-compute-fractional-powers}
\varLambda(t)^{-\alpha}=\dfrac{\sin\pi\alpha}{\pi}\int_0^\infty\lambda^{-\alpha}(\lambda I+\varLambda(t))^{-1}d\lambda.
\end{equation}

\noindent
v) For each $\alpha\in(0,1)$, the operator $\varLambda(t)^\alpha$ is a negative generator of an analytic $C^0$-semigroup $\{e^{-\varLambda(t)^\alpha \tau}:\tau\geqslant 0\}$.

\noindent
vi) Given any $0<\alpha<1$ we have 
\begin{equation*}
\varLambda(t)^{-\alpha}=
\begin{bmatrix}
\mu(t)^{-\frac{\alpha}{2}}\cos\frac{\pi\alpha}{2}A^{-\frac{\alpha}{2}} & \mu(t)^{\frac{-1-\alpha}{2}}\sin\frac{\pi\alpha}{2}A^{\frac{-1-\alpha}{2}} \\ \\
-\mu(t)^{\frac{1-\alpha}{2}} \sin\frac{\pi\alpha}{2}A^{\frac{1-\alpha}{2}} & \mu(t)^{-\frac{\alpha}{2}}\cos\frac{\pi\alpha}{2}A^{-\frac{\alpha}{2}}
\end{bmatrix}
\end{equation*}
and
\begin{equation}\label{pot-a}
\varLambda(t)^\alpha= \begin{bmatrix}
\mu(t)^{\frac{\alpha}{2}} \cos\frac{\pi\alpha}{2}A^{\frac{\alpha}{2}} & -\mu(t)^{\frac{-1+\alpha}{2}}\sin\frac{\pi\alpha}{2}A^{\frac{-1+\alpha}{2}} \\ \\
\mu(t)^{\frac{1+\alpha}{2}}\sin\frac{\pi\alpha}{2}A^{\frac{1+\alpha}{2}} & \mu(t)^{\frac{\alpha}{2}} \cos\frac{\pi\alpha}{2}A^{\frac{\alpha}{2}}
\end{bmatrix},
\end{equation}
with domain $D(\varLambda(t)^\alpha)=X^{\frac{1+\alpha}{2}}\times X^{\frac{\alpha}{2}}$.

\noindent
vii) For each $\alpha\in(0, 1]$ and for every $t\in\mathbb{R}$, the spectrum of $-\Lambda(t)^{\alpha}$ is a point spectrum consisting of eigenvalues 
\begin{equation*}%\label{eigenvalues}
\lambda^{\pm}_{\alpha,n}(t)=e^{\pm i\frac{\pi(2-\alpha)}{2}}\mu(t)^{\frac{\alpha}{2}}\nu^{\frac{\alpha}{2}}_{n}, \quad n\in\mathbb{N},
\end{equation*}
where $\{\nu_n\}_{n\in\mathbb{N}}$ denotes the ordered sequence of eigenvalues of $A$ including their multiplicity.

\noindent
viii) The family $\varLambda(t)^{-\alpha}$ converges to $\varLambda(t)^{-1}$ in $\mathscr{L}(Y_t)$ as $\alpha\nearrow 1$.

\noindent
ix) For each $\left[\begin{smallmatrix} \varphi \\ \psi \end{smallmatrix}\right]\in Y_t^1$,
$$
\varLambda(t)^\alpha \left[\begin{smallmatrix} \varphi \\ \psi \end{smallmatrix}\right] \stackrel{}{\to} \varLambda(t) \left[\begin{smallmatrix} \varphi \\ \psi \end{smallmatrix}\right] \ \text{ in } \ Y_t \ \text{ as } \ \alpha\nearrow1.
$$

\end{lemma}

\proof
Parts $i)$ and $ii)$ are consequence of (\ref{eq:A}) and (\ref{eq:Lambda}). Part $iii)$ comes froms $ii)$ and Stone's theorem (see Pazy \cite[Theorem 10.8]{Pa}). Parts $iv)$ and $v)$ follows from Kato \cite[Theorems~1,~2]{K1}. Concerning part $vi)$ note that given $\lambda\in\mathbb{C}$ we have
\begin{equation*}
\lambda I+\varLambda(t)=\begin{bmatrix} \lambda I & -I\\
\mu(t)A & \lambda I
\end{bmatrix}
\end{equation*}
and
\begin{equation*}\label{resolvente-A}
(\lambda I+\varLambda(t))^{-1}= \begin{bmatrix} \lambda(\lambda^2 I+\mu(t)A)^{-1} & (\lambda^2 I+\mu(t)A)^{-1}  \\ \\
-\mu(t)A(\lambda^2 I+\mu(t)A)^{-1}& \lambda(\lambda^2 I+\mu(t)A)^{-1}
\end{bmatrix} \ \text{ for all $\lambda\in\rho(-\varLambda(t))$}.
\end{equation*}
Due to (\ref{eq:to-compute-fractional-powers}) for any $0<\alpha<1$, we get
\begin{equation*}
\varLambda(t)^{-\alpha}=
\begin{bmatrix}
\mu(t)^{-\frac{\alpha}{2}}\cos\frac{\pi\alpha}{2}A^{-\frac{\alpha}{2}} & \mu(t)^{\frac{-1-\alpha}{2}}\sin\frac{\pi\alpha}{2}A^{\frac{-1-\alpha}{2}} \\ \\
-\mu(t)^{\frac{1-\alpha}{2}} \sin\frac{\pi\alpha}{2}A^{\frac{1-\alpha}{2}} & \mu(t)^{-\frac{\alpha}{2}}\cos\frac{\pi\alpha}{2}A^{-\frac{\alpha}{2}}
\end{bmatrix}
\end{equation*}
which leads to (\ref{pot-a}).

Part $vii)$ was proved in \cite[Theorem 4.8]{BN}. Part $viii)$ follows from Amann \cite[Theorem III.4.6.2]{Am}. Finally, to prove part $ix)$ we fix $\left[\begin{smallmatrix} \varphi \\ \psi \end{smallmatrix}\right]\in Y_t^1$ and observe that
$$
\cos\dfrac{\pi\alpha}{2}A^{\frac{\alpha}{2}} \varphi \stackrel{X^\frac12}\to 0 \quad \text{ and } \quad \cos\dfrac{\pi\alpha}{2}A^{\frac{\alpha}{2}} \psi \stackrel{X}\to 0 \ \text{ as } \alpha \nearrow 1,
$$
because on the one hand $\cos\dfrac{\pi\alpha}{2}\to 0$ and, on the other, due to $viii)$,
$$A^\frac12 A^{\frac{\alpha}{2}} \varphi = A^{-\frac{1-\alpha}2}A^1 \varphi\stackrel{X}\to A^1 \varphi \quad \text{ and } \quad  A^{\frac{\alpha}{2}} \psi = A^{-\frac{1-\alpha}2}A^\frac12 \psi\stackrel{X}\to A^\frac12 \varphi \ \text{ as } \ \alpha\nearrow 1.$$
In a similar manner, using that $\sin\dfrac{\pi\alpha}{2}\to 1$ as $\alpha\nearrow 1$ and that, due to $viii)$,
$$
A^{\frac{1+\alpha}{2}} \varphi  = A^{\frac{-1+\alpha}{2}} A^1 \varphi \stackrel{X}\to A\phi \quad  \text{ and } \quad
A^{\frac{-1+\alpha}{2}}A^\frac12\psi \stackrel{X}\to A^\frac12\psi,
$$
we get
$$
\sin\dfrac{\pi\alpha}{2}A^{\frac{1+\alpha}{2}} \varphi \stackrel{X}\to A\phi \quad \text{ and } \ \quad
\sin\dfrac{\pi\alpha}{2}A^{\frac{-1+\alpha}{2}}\psi \stackrel{X^\frac12}\to \psi \ \text{ as } \alpha \nearrow 1.
$$
Hence
$$
\begin{bmatrix}
\mu(t)^{\frac{\alpha}{2}} \cos\frac{\pi\alpha}{2}A^{\frac{\alpha}{2}} & -\mu(t)^{\frac{-1+\alpha}{2}}\sin\frac{\pi\alpha}{2}A^{\frac{-1+\alpha}{2}} \\ \\
\mu(t)^{\frac{1+\alpha}{2}}\sin\frac{\pi\alpha}{2}A^{\frac{1+\alpha}{2}} & \mu(t)^{\frac{\alpha}{2}} \cos\frac{\pi\alpha}{2}A^{\frac{\alpha}{2}}
\end{bmatrix} \begin{bmatrix} \varphi \\ \psi \end{bmatrix} \stackrel{Y_t}\to \begin{bmatrix} 0 & -I\\ \mu(t)A & 0
\end{bmatrix}  \begin{bmatrix} \phi
\\ \varphi \end{bmatrix}
\ \text{for each $ \begin{bmatrix} \varphi \\ \psi \end{bmatrix} \in Y_t^1$},
$$
which gives the result.
\qed

\medskip

\begin{remark}
Thanks to Bezerra and Nascimento \cite[Remark 3.1]{BN}, for any $0<\alpha<1$ we have the following identity  
\[
(\mu(t)A)^\alpha=\mu(t)^\alpha A^\alpha,\ \mbox{for all}\ t\in\mathbb{R}.
\]
\end{remark}

\medskip

Next, we present some results of the relationship between our fractional powers and interpolation spaces.

\begin{lemma}
The extrapolation space of $X^\frac12\times X$ generated by $\varLambda(t_0)$ coincides with $X\times X^{-\frac12}$, for any $t_0\in\mathbb{R}$.
\end{lemma}

\proof
Let $t_0\in\mathbb{R}$, using Lemma~\ref{lem:useful:properties} $i)$, we have
$$
\bigl\|\varLambda(t_0)^{-1} \left[\begin{smallmatrix}\varphi\\ \psi \end{smallmatrix}\right] \bigr\|_{X^\frac12\times X}= \Bigl\| \left[\begin{smallmatrix}  \mu(t_0)^{-1}A^{-1} \psi\\ -\varphi
\end{smallmatrix}\right] \Bigr\|_{X^\frac12\times X} = \bigl\| \left[\begin{smallmatrix}\varphi\\ \psi \end{smallmatrix}\right] \bigr\|_{X\times X^{-\frac12}} \quad \text{ for any } \ \left[\begin{smallmatrix}\varphi\\ \psi \end{smallmatrix}\right] \in X^\frac12\times X
$$
and taking the completion we get the result (see Carvalho and Cholewa \cite[Lemma~2]{CC}).
\qed

\begin{lemma}
Let $\alpha\in(0,1]$ be fixed and let $\{E^\theta_t(\alpha), \theta\in[-1,\infty)\}$ be the extrapolated fractional power scale of order $1$ generated by $(X^\frac12\times X, \varLambda(t)^\alpha)$. Then
\begin{equation}\label{eq:characterization}
E^\theta_t(\alpha)= X^\frac{1+\alpha\theta}2\times X^\frac{\alpha\theta}2 \ \text{ for each } \ \theta\in[-1,1]\text{ and }t\in\mathbb{R}.
\end{equation}
\end{lemma}
\proof The proof is similar to  Bezerra et al. \cite[Lemma 2.3]{BCCN}.  \qed

\section{The fractional oscillon equations}\label{Sec:frac-osc-eq}

In this section, we obtain a model of oscillon equation via fractional powers of order $\alpha\in(0,1)$. 

\begin{definition}
Given $\left[\begin{smallmatrix} u_\tau\\ v_\tau\end{smallmatrix}\right]\in X\times X^{-\frac12}$ we say that $\left[\begin{smallmatrix} u\\ v\end{smallmatrix}\right]$ is a global mild solution of \eqref{yyy} provided that $\left[\begin{smallmatrix} u\\ v\end{smallmatrix}\right]\in C([\tau,\infty), X\times X^{-\frac12})$,
$f^e(u)\in C([\tau,\infty), X^{-\frac12})$ and $\left[\begin{smallmatrix} u\\ v\end{smallmatrix}\right]$ satisfies for $t>\tau$ the integral equation
\begin{equation*}\label{eq:vocf}
\left[\begin{smallmatrix} u(t)\\ v(t)\end{smallmatrix}\right]=U_1(t,\tau)\left[\begin{smallmatrix} u_\tau\\ v_\tau\end{smallmatrix}\right]
+\int_\tau^t U_1(t,s)\left[\begin{smallmatrix} 0\\ f^e(u(s))-\omega(s)v(s)\end{smallmatrix}\right]ds,
\end{equation*}
where
\[
U_1(t,\tau) = e^{-(t-\tau)\varLambda(\tau)} + \int_{\tau}^{t}U_1(t,s)[\varLambda(\tau) -\varLambda(s)]e^{-(s-\tau)\varLambda(\tau)}ds,\quad t>\tau.
\]
\end{definition}

Initially, we obtain a result on local well posedness of the Cauchy problem  \eqref{ffhghgf} in spaces $X^\frac{\alpha}2\times X^{\frac{\alpha-1}2}$ for some $0<\alpha<1$. Moreover, we show some dissipativity property of those solutions.

\begin{theorem}\label{th:local-well-posedness-of-perturbed-problem}
Suppose that \eqref{growth-condition} and \eqref{eq:supercritical} hold, $\omega$ is $(\zeta, \kappa_0)-$H\"{o}lder continuous in $\mathbb{R}$ satisfying \eqref{132r} and fix any number $s$ satisfying $\frac{N}2(1-\frac1\rho) -\frac1\rho<s<1$. Then, for each
$\alpha\in[\frac{N}2(\rho-1) - \rho s,1)$ the following hold:
\begin{itemize}
\item[$i)$] For any $\left[\begin{smallmatrix}u_{\tau\alpha}\\ v_{\tau\alpha} \end{smallmatrix}\right]\in X^\frac{s}2\times X^\frac{s-1}2$ there exists a unique mild solution 
\[
\left[\begin{smallmatrix}u^{\alpha}\\ v^{\alpha} \end{smallmatrix}\right]\in C([\tau,\xi_{u_{\tau\alpha}, v_{\tau\alpha}}), X^\frac{s}2\times X^{\frac{s-1}2})
\]
of \eqref{ffhghgf} defined on a maximal interval of existence $[\tau,\xi_{u_{\tau\alpha}, v_{\tau\alpha}})$. This solution depends continuously on the initial data and satisfies a blow up alternative in $X^\frac{s}2\times X^\frac{s-1}2$. In particular, if $\|\left[\begin{smallmatrix}u^{\alpha}\\ v^{\alpha} \end{smallmatrix}\right]\|_{X^\frac{s}2\times X^\frac{s-1}2}$-norm remains bounded as long as the solution exists then $\xi_{u_{\tau\alpha}, v_{\tau\alpha}}=\infty$.
\item[$ii)$] The solution in part $i)$ above is a regular solution. Namely,
\begin{equation}\label{eq:regularity-of-solution}
\left[\begin{smallmatrix}u^{\alpha}\\ v^{\alpha} \end{smallmatrix}\right]\in C((\tau,\xi_{u_{\tau\alpha}, v_{\tau\alpha}}), X^\frac{1+\alpha}2\times X^\frac\alpha2)\cap C^1((\tau,\xi_{u_{\tau\alpha}, v_{\tau\alpha}}), X^\frac{1+\sigma}2\times X^\frac\sigma2) \ \text{ for each $\sigma<\alpha$}
\end{equation}
and $\left[\begin{smallmatrix}u^{\alpha}\\ v^{\alpha} \end{smallmatrix}\right]$ satisfies \eqref{ffhghgf}.
\item[$iii)$] Actually, for any set $B$ bounded in $X^\frac{s}2\times X^\frac{s-1}2$ there is a certain time $\tau_B>\tau\in\mathbb{R}$ such that for each $\left[\begin{smallmatrix}u_{\tau\alpha}\\ v_{\tau\alpha} \end{smallmatrix}\right]\in B$ the solution $\left[\begin{smallmatrix}u^{\alpha}\\ v^{\alpha} \end{smallmatrix}\right]$ through $\left[\begin{smallmatrix}u_{\tau\alpha}\\ v_{\tau\alpha} \end{smallmatrix}\right]$ in part $i)$ exists (at least) until $\tau_B$ and given any $\xi\in(\tau,\tau_B]$ there is a positive constant $M=M(\xi,B)$ such that
\begin{equation}\label{eq:uniform-regularization}
\sup_{\left[\begin{smallmatrix}u_{\tau\alpha}\\ v_{\tau\alpha} \end{smallmatrix}\right]\in B} \Bigl\| \left[\begin{smallmatrix}u^{\alpha}(\xi)\\ v^{\alpha} (\xi) \end{smallmatrix}\right]\Bigr\|_{X^\frac{1+\alpha}2\times X^\frac\alpha2} \leqslant  M.
\end{equation}
\end{itemize}
\end{theorem}

\proof  Part $i)$ follows from Theorem \ref{theo:10r} $i)$, Lemma \ref{lemma:CD02} and Caraballo et al. \cite[Theorem 2.3]{CCLR}. To prove $ii)$ we refer to Bezerra et al. \cite[Theorem 1.2 ({\rm ii)}]{BCCN}. Finally, item $iii)$ follows from Carvalho and Cholewa \cite[Theorem 5]{CC1}.
\qed

\medskip

The fractional abstract model \eqref{ffhghgf} can be seen as we will show below as a Cauchy problem associated with a space-fractional PDE.

Thanks to \eqref{eq:Lambdaf} and \eqref{pot-a}, we rewrite the Cauchy problem \eqref{ffhghgf} as
\begin{equation}\label{6f}
\begin{cases}
v^\alpha+\mu(t)^{\frac{\alpha}{2}} \cos\frac{\pi\alpha}{2}A^{\frac{\alpha}{2}}u^\alpha-\mu(t)^{\frac{-1+\alpha}{2}}\sin\frac{\pi\alpha}{2}A^{\frac{-1+\alpha}{2}}v^\alpha=0,& t>\tau,\\
v_t^\alpha+\mu(t)^{\frac{1+\alpha}{2}}\sin\frac{\pi\alpha}{2}A^{\frac{1+\alpha}{2}} u^\alpha+ \mu(t)^{\frac{\alpha}{2}} \cos\frac{\pi\alpha}{2}A^{\frac{\alpha}{2}}v^\alpha=f(u^\alpha)-\omega(t)v^\alpha,& t>\tau,\\
u^\alpha(\tau)=u_{\tau\alpha},\\
v^\alpha(\tau)=v_{\tau\alpha},\\
\end{cases}
\end{equation}
where $v^\alpha=u_{t}^\alpha$.

We can see from the first equation in \eqref{6f} that
\begin{equation}\label{6fb}
\mu(t)^{\frac{\alpha}{2}}\sin\frac{\pi\alpha}{2}\cos\frac{\pi\alpha}{2}A^{\frac{\alpha}{2}}v^\alpha-\mu(t)^{\frac{1+\alpha}{2}} \cos^2\frac{\pi\alpha}{2}A^{\frac{1+\alpha}{2}}u^\alpha=\mu(t)^{\frac{1}{2}}\cos\frac{\pi\alpha}{2} A^{\frac{1}{2}}v^\alpha,
\end{equation}
and 
\[
\sin\frac{\pi\alpha}{2}v^\alpha=A^{\frac{1-\alpha}{2}}[\mu(t)^{\frac{1-\alpha}{2}}v^\alpha+\mu(t)^{\frac{1}{2}} \cos\frac{\pi\alpha}{2}A^{\frac{\alpha}{2}}u^\alpha],
\]
which implies 
\begin{equation}\label{7fb}
\begin{split}
\sin\frac{\pi\alpha}{2}v_t^\alpha&=A^{\frac{1-\alpha}{2}}\Big[\Big(\frac{1-\alpha}{2}\Big)\mu(t)^{\frac{-1-\alpha}{2}}\mu'(t)v^\alpha+\mu(t)^{\frac{1-\alpha}{2}}v_t^\alpha\Big]\\
&+A^{\frac{1-\alpha}{2}}\Big[\frac{1}{2}\mu(t)^{-\frac{1}{2}}\mu'(t)\cos\frac{\pi\alpha}{2} A^{\frac{\alpha}{2}}u^\alpha+\mu(t)^{\frac{1}{2}}\cos\frac{\pi\alpha}{2} A^{\frac{\alpha}{2}}v^\alpha\Big].
\end{split}
\end{equation}
From the second equation in \eqref{6f}, we obtain
\begin{equation}\label{8fb}
\begin{split}
&\sin\frac{\pi\alpha}{2}f(u^\alpha)\\
&=\sin\frac{\pi\alpha}{2}\omega(t)v^\alpha+\sin\frac{\pi\alpha}{2}v_t^\alpha+\mu(t)^{\frac{1+\alpha}{2}}\sin^2\frac{\pi\alpha}{2}A^{\frac{1+\alpha}{2}} u^\alpha+ \mu(t)^{\frac{\alpha}{2}} \sin\frac{\pi\alpha}{2}\cos\frac{\pi\alpha}{2}A^{\frac{\alpha}{2}}v^\alpha\\
&=\sin\frac{\pi\alpha}{2}\omega(t)v^\alpha+\sin\frac{\pi\alpha}{2}v^{\alpha}_t+\mu(t)^{\frac{1+\alpha}{2}} A^{\frac{1+\alpha}{2}} u^\alpha+\mu(t)^{\frac{\alpha}{2}} \sin\frac{\pi\alpha}{2}\cos\frac{\pi\alpha}{2}A^{\frac{\alpha}{2}}v^\alpha\\
&-\mu(t)^{\frac{1+\alpha}{2}}\cos^2\frac{\pi\alpha}{2}A^{\frac{1+\alpha}{2}} u^\alpha \\
\end{split}
\end{equation}
and by \eqref{6fb} and \eqref{8fb}, we have
\begin{equation}\label{9fb}
\sin\frac{\pi\alpha}{2}f(u^\alpha)=\sin\frac{\pi\alpha}{2}\omega(t)v^\alpha+\sin\frac{\pi\alpha}{2}v_t^\alpha+\mu(t)^{\frac{1+\alpha}{2}} A^{\frac{1+\alpha}{2}} u^\alpha+\mu(t)^{\frac{1}{2}}\cos\frac{\pi\alpha}{2} A^{\frac{1}{2}}v^\alpha.
\end{equation}

Thus, using \eqref{7fb} and \eqref{9fb}, we conclude that
\begin{equation}\label{10fb}
\begin{split}
&\mu(t)^{\frac{1-\alpha}{2}}A^{\frac{1-\alpha}{2}}u_{tt}^\alpha+\frac{1}{2}\mu(t)^{-\frac{1}{2}}\mu'(t)\cos\frac{\pi\alpha}{2} A^{\frac{1}{2}} u^\alpha+\mu(t)^{\frac{1+\alpha}{2}} A^{\frac{1+\alpha}{2}} u^\alpha+2\mu(t)^{\frac{1}{2}}\cos\frac{\pi\alpha}{2} A^{\frac{1}{2}}u_t^\alpha\\
&+\Big(\frac{1-\alpha}{2}\Big)\mu(t)^{\frac{-1-\alpha}{2}}\mu'(t)A^{\frac{1-\alpha}{2}}u_t^\alpha +\sin\frac{\pi\alpha}{2}\omega(t)u_t^\alpha=\sin\frac{\pi\alpha}{2}f(u^\alpha).
\end{split}
\end{equation}

\begin{remark}
For $n=1$, $\alpha=1$, $\omega\equiv H>0$ and $\mu(t)=e^{-2Ht}$ the fractional PDE in \eqref{10fb} it becomes the oscillon equation in \eqref{asdar}, which was treated in Di Plino, Duane and Temam \cite{DDT1} with periodic boundary conditions. 

For $n=3$ and $\alpha=1$ the fractional PDE in \eqref{10fb} it becomes the oscillon equation in \eqref{asdar01}, which was treated in Di Plino, Duane and Temam \cite{DDT2}. 
\end{remark}

\begin{remark}
Since \eqref{10fb} rewrites  for $\mathbb{A}=A^\alpha$ as 
\begin{equation}\label{12fbr}
\begin{split}
\mu(t)^{\frac{1-\alpha}{2}}u_{tt}^\alpha&+\frac{1}{2}\mu(t)^{-\frac{1}{2}}\mu'(t)\cos\frac{\pi\alpha}{2} \mathbb{A}^{\frac{1}{2}}u^\alpha +\mu(t)^{\frac{1+\alpha}{2}} \mathbb{A} u^\alpha\\
&+\Big(\frac{1-\alpha}{2}\Big)\mu(t)^{\frac{-1-\alpha}{2}}\mu'(t) u_t^\alpha+2\mu(t)^{\frac{1}{2}}\cos\frac{\pi\alpha}{2} \mathbb{A}^{\frac{1}{2}}u_t^\alpha+\sin\frac{\pi\alpha}{2}\omega(t)\mathbb{A}^{-\frac{1-\alpha}{2\alpha}}u_t^\alpha\\
&=\sin\frac{\pi\alpha}{2}\mathbb{A}^{-\frac{1-\alpha}{2\alpha}} f(u^\alpha),
\end{split}
\end{equation}
this latter equation can be viewed as a parabolic-type approximation of \eqref{xxx} (see  Carvalho and Cholewa \cite{CC,CC1} and Chen and Triggiani \cite{CSTR} for the extensive studies of the strongly damped wave equations).

For $\mu\equiv1$ and   $\omega\equiv a$ with $a>0$, the fractional PDE in \eqref{12fbr} it becomes the autonomous fractional PDE which was treated in Bezerra et al. \cite[Remark 3.5]{BCCN}.
\end{remark}

\subsection{Lyapunov functionals associated with perturbed problems}

 The aim of this section is define a Lyapunov functional associated with our perturbed problem \eqref{ffhghgf}. With this purpose, following Hale \cite{hale}, it follows from \eqref{CondDissip} that, for each $\varepsilon > 0$, there is a constant $C_{\varepsilon}>0$ such that

\begin{equation}\label{1558r}
\langle f(u), u \rangle\leqslant \varepsilon\|u\|^{2}+C_{\varepsilon}\end{equation}
and 
\begin{equation}\label{1559r}
|V(u)|\leqslant \varepsilon\|u\|^{2}+C_{\varepsilon},
\end{equation}
where $V(u):=\int_\Omega\int_0^{u}f(\sigma)d\sigma dx$, for each $u\in X$ such that $f(u)u\in L^1(\Omega)$ and $\int_0^{s}f(\sigma)d\sigma\in L^1(\Omega)$.

\medskip
We will write problem \eqref{6f} in our abstract framework. For $0<\alpha\leqslant 1$ and $t\in\mathbb{R}$, we introduce the Banach space
\[
\mathcal{H}^\alpha_t:=X^{\frac{1+\alpha}{4}}\times X^{\frac{1-\alpha}{4}}
\]
with norm
\begin{equation}\label{3.8}
\|\left[\begin{smallmatrix} u^{\alpha} \\ v^{\alpha} \end{smallmatrix}\right]\|_{\mathcal{H}^\alpha_t}=\mu(t)^{\frac{1+\alpha}{4}}\|u^{\alpha}\|_{X^{\frac{1+\alpha}{4}}}+\|u^{\alpha}\|_{X}+\mu(t)^{\frac{1-\alpha}{4}}\|v^{\alpha}\|_{X^{\frac{1-\alpha}{4}}},\quad \forall (u^{\alpha},v^{\alpha})\in \mathcal{H}^\alpha_t.
\end{equation}
For simplicity, we set $\mathcal{H}^\alpha=\mathcal{H}^\alpha_0$.

For some of the proofs below, it will be convenient to use the natural energy of the problem at time $t$
\[
\mathscr{E}_{\mathcal{H}^\alpha_t}(u^{\alpha},v^{\alpha})=\mu(t)^{\frac{1+\alpha}{2}}\|u^\alpha\|^2_{X^{\frac{1+\alpha}{4}}}+\|u^\alpha\|^2+\mu(t)^{\frac{1-\alpha}{2}}\|u_t^\alpha\|^2_{X^{\frac{1-\alpha}{4}}}
\]
 instead of $\mathcal{H}^\alpha_t$-norm.  Note that
\[
\mathscr{E}_{\mathcal{H}^\alpha_t}(u^{\alpha},v^{\alpha})\leqslant \|\left[\begin{smallmatrix} u^{\alpha} \\ v^{\alpha} \end{smallmatrix}\right]\|_{\mathcal{H}^\alpha_t}^2\quad\mbox{and}\quad \|\left[\begin{smallmatrix} u^{\alpha} \\ v^{\alpha} \end{smallmatrix}\right]\|_{\mathcal{H}^\alpha_t}^2\leqslant C\mathscr{E}_{\mathcal{H}^\alpha_t}(u^{\alpha},v^{\alpha}),
\]
for some constant $C>0$, and then the energy $\mathscr{E}_{\mathcal{H}^\alpha_t}(\cdot)$ is equivalent to the norm $\|\cdot\|_{\mathcal{H}^\alpha_t}$.

\begin{definition}
A family $\mathcal{B}=\{\mathcal{B}(t)\subset \mathcal{H}^\alpha_t: t\in\mathbb{R}\}$ is pullback-bounded if and only if 
\[
\sup_{s\in(-\infty,t]}\sup_{z\in\mathcal{B}(t)}\mathscr{E}_{\mathcal{H}^\alpha_t}(z)<\infty,\quad\forall t\in\mathbb{R};
\]
\end{definition}

\begin{definition}
A sequence $\{z_n\}\subset \mathcal{H}^\alpha_t$ converges to $z\in \mathcal{H}^\alpha_t$ if and only if $\mathscr{E}_{\mathcal{H}^\alpha_t}(z_n-z)$ converges to zero as $n\to\infty$.
\end{definition}

Next, we obtain a suitable a-priori estimate for the solution of \eqref{10fb}.

\begin{lemma}\label{inequality-geral} 
Let $s\in(\frac{N}2(1-\frac1\rho) -\frac1\rho,1)$, $\left[\begin{smallmatrix}u_{\tau\alpha}\\ v_{\tau\alpha} \end{smallmatrix}\right]\in X^\frac{s}2\times X^\frac{s-1}2$ and $S_{\alpha}(t,\tau)\left[\begin{smallmatrix}u_{\tau\alpha}\\ v_{\tau\alpha} \end{smallmatrix}\right]=\left[\begin{smallmatrix}u^{\alpha}\\ v^{\alpha} \end{smallmatrix}\right]$ be the solution of \eqref{10fb} with initial time $\tau\in\mathbb{R}$ and  initial data $\left[\begin{smallmatrix}u_{\tau\alpha}\\ v_{\tau\alpha} \end{smallmatrix}\right]$, for all $\alpha\in[\frac{N}2(\rho-1) - \rho s,1)$.

The following a-priori estimate holds:
\begin{equation}\label{160r}
\mathscr{E}_{\mathcal{H}^\alpha_t}(S_{\alpha}(t,\tau)\left[\begin{smallmatrix}u_{\tau\alpha}\\ v_{\tau\alpha} \end{smallmatrix}\right])\leqslant M_0 \mathscr{E}_{\mathcal{H}^\alpha_\tau}(\left[\begin{smallmatrix}u_{\tau\alpha}\\ v_{\tau\alpha} \end{smallmatrix}\right])e^{-\varepsilon_{\omega}(t)(t-\tau)}+M_1, \quad\forall t\geqslant \tau,
\end{equation}
where $M_0$ and $M_1$ are positive constants.
\end{lemma}

\proof We denote by $(u^{\alpha}(t),u^{\alpha}_t(t))$, the solution of \eqref{10fb} with initial time $\tau\in\mathbb{R}$ and initial condition  $\left[\begin{smallmatrix}u_{\tau\alpha}\\ v_{\tau\alpha} \end{smallmatrix}\right]\in X^\frac{s}2\times X^\frac{s-1}2$, which we assume to be sufficiently regular (see \eqref{eq:regularity-of-solution}).

Firstly, multiplying  the  equation in \eqref{10fb}  by $u^{\alpha}_t$ and integrating over $\Omega$, we have 
\[
\begin{split}
\mu(t)^{\frac{1-\alpha}{2}}&\langle A^{\frac{1-\alpha}{2}}u_{tt}^\alpha,u_t^\alpha\rangle+\Big\langle\Big[\frac{1}{2}\mu(t)^{-\frac{1}{2}}\mu'(t)\cos\frac{\pi\alpha}{2} A^{\frac{1}{2}} +\mu(t)^{\frac{1+\alpha}{2}} A^{\frac{1+\alpha}{2}} \Big]u^\alpha,u_t^\alpha\Big\rangle\\
&+\Big\langle\Big[\Big(\frac{1-\alpha}{2}\Big)\mu(t)^{\frac{-1-\alpha}{2}}\mu'(t)A^{\frac{1-\alpha}{2}} +2\mu(t)^{\frac{1}{2}}\cos\frac{\pi\alpha}{2} A^{\frac{1}{2}}+\sin\frac{\pi\alpha}{2}\omega(t)\Big]u_t^\alpha,u_t^\alpha\Big\rangle\\
&=\sin\frac{\pi\alpha}{2}\langle f(u^\alpha),u_t^\alpha\rangle.
\end{split}
\]
Therefore
\begin{equation}\label{12fb}
\begin{split}
\dfrac{d}{dt}&\Big[\mu(t)^{\frac{1+\alpha}{2}} \|u^\alpha\|^2_{X^{\frac{1+\alpha}{4}} }+\mu(t)^{\frac{1-\alpha}{2}}\|u_t^\alpha\|^2_{X^{\frac{1-\alpha}{4}}}-2\sin\frac{\pi\alpha}{2}\int_\Omega\int_0^{u^\alpha}f(s)dsdx\Big]\\
&-\Big(\frac{1+\alpha}{2}\Big)\mu(t)^{\frac{-1-\alpha}{2}}\mu'(t)\|u_t^\alpha\|^2_{X^{\frac{1-\alpha}{4}}}-\Big(\frac{1+\alpha}{2}\Big)\mu(t)^{\frac{-1+\alpha}{2}}\mu'(t)\|u^\alpha\|^2_{X^{\frac{1+\alpha}{4}}}\\
&+4\mu(t)^{\frac{1}{2}}\cos\frac{\pi\alpha}{2}\|u_t^\alpha\|^2_{X^{\frac{1}{4}}}+\sin\frac{\pi\alpha}{2}(\omega(t)+\mu(t)^{-1}\mu'(t))\|u_t^\alpha\|^2=0.
\end{split}
\end{equation}
On the other hand, multiplying  the  equation in \eqref{10fb}  by $u$ and integrating over $\Omega$, we get
\[
\begin{split}
\mu&(t)^{\frac{1-\alpha}{2}}\langle A^{\frac{1-\alpha}{2}}u_{tt}^\alpha,u^\alpha\rangle+\Big\langle\Big[\frac{1}{2}\mu(t)^{-\frac{1}{2}}\mu'(t)\cos\frac{\pi\alpha}{2} A^{\frac{1}{2}} +\mu(t)^{\frac{1+\alpha}{2}} A^{\frac{1+\alpha}{2}} \Big]u^\alpha,u^\alpha\Big\rangle\\
&+\Big\langle\Big[\Big(\frac{1-\alpha}{2}\Big)\mu(t)^{\frac{-1-\alpha}{2}}\mu'(t)A^{\frac{1-\alpha}{2}} +2\mu(t)^{\frac{1}{2}}\cos\frac{\pi\alpha}{2} A^{\frac{1}{2}}+\sin\frac{\pi\alpha}{2}\omega(t)\Big]u_t^\alpha,u^\alpha\Big\rangle\\
&=\sin\frac{\pi\alpha}{2}\langle f(u^\alpha),u^\alpha\rangle.
\end{split}
\]
Hence, for $\varepsilon> 0$ to be determined later and using \eqref{1558r}, we obtain
\begin{equation}\label{2r}
\begin{split}
\dfrac{d}{dt}\Big[&2\mu(t)^{\frac{1}{2}}\cos\frac{\pi\alpha}{2} \|u^\alpha\|^2_{X^{\frac{1}{4}} }+\sin\frac{\pi\alpha}{2}\omega(t)\|u^\alpha\|^2+2\mu(t)^{\frac{1-\alpha}{2}}\langle A^{\frac{1-\alpha}{4}}u_t^{\alpha}, A^{\frac{1-\alpha}{4}} u^{\alpha}\rangle\Big]\\
&+2\mu(t)^{\frac{1+\alpha}{2}}\|u^\alpha\|^2_{X^{\frac{1+\alpha}{4}}} -2\mu(t)^{\frac{1-\alpha}{2}}\|u_t^\alpha\|^2_{X^{\frac{1-\alpha}{4}}}-\sin\frac{\pi\alpha}{2}\omega^{\prime}(t)\|u^\alpha\|^2\\
&=2\sin\frac{\pi\alpha}{2}\langle f(u^{\alpha}), u^{\alpha}\rangle\leqslant 2\varepsilon\sin\frac{\pi\alpha}{2} \|u^\alpha\|^{2} + 2 C_{\varepsilon}\sin\frac{\pi\alpha}{2}.
\end{split}
\end{equation}
We add \eqref{12fb} to $2\varepsilon-$times \eqref{2r}. Thus,
\begin{equation}\label{10r}
\dfrac{d}{dt}\Phi_{\alpha}(u^\alpha,u^\alpha_t)+2\varepsilon\Phi_{\alpha}(u^\alpha,u^\alpha_t)+\Phi_{*,\alpha}(u^\alpha,u^\alpha_t)+\sin\frac{\pi\alpha}{2}\omega(t)\|u^{\alpha}_t\|^2\leqslant 4\varepsilon C_{\varepsilon},
\end{equation}
where
\begin{equation}\label{funcional01}
\begin{split}
& \Phi_{\alpha}(u^\alpha,u^\alpha_t)\\
&=\mu(t)^{\frac{1+\alpha}{2}}\|u^\alpha\|^2_{X^{\frac{1+\alpha}{4}}}+\mu(t)^{\frac{1-\alpha}{2}}\|u_t^\alpha\|^2_{X^{\frac{1-\alpha}{4}}}-2\sin\frac{\pi\alpha}{2}V(u^{\alpha})\\
&+2\varepsilon\Big[2\mu(t)^{\frac{1}{2}}\cos\frac{\pi\alpha}{2} \|u^\alpha\|^2_{X^{\frac{1}{4}} }+\sin\frac{\pi\alpha}{2}\omega(t)\|u^\alpha\|^2+2\mu(t)^{\frac{1-\alpha}{2}}\langle A^{\frac{1-\alpha}{4}}u_t^{\alpha}, A^{\frac{1-\alpha}{4}} u^{\alpha}\rangle\Big]
\end{split}
\end{equation}
and
\[
\begin{split}
&\Phi_{*,\alpha}(u^\alpha,u^\alpha_t)\\
&=\Big[-\Big(\frac{1+\alpha}{2}\Big)\mu(t)^{\frac{-1-\alpha}{2}}\mu'(t)-6\varepsilon\mu(t)^{\frac{1-\alpha}{2}}\Big] \|u^{\alpha}_t\|^{2}_{X^{\frac{1-\alpha}{4}}}+4\mu(t)^{\frac{1}{2}}\cos\frac{\pi\alpha}{2}\|u^{\alpha}_t\|^{2}_{X^{\frac{1}{4}}}\nonumber\\
&-8\varepsilon^{2}\mu(t)^{\frac{1}{2}}\cos\frac{\pi\alpha}{2}\|u^{\alpha}\|^{2}_{X^{\frac{1}{4}}}+\Big[-\Big(\frac{1+\alpha}{2}\Big)\mu(t)^{\frac{-1+\alpha}{2}}\mu'(t)+2\varepsilon\mu(t)^{\frac{1+\alpha}{2}}\Big] \|u^{\alpha}\|^{2}_{X^{\frac{1+\alpha}{4}}}\nonumber\\
&+\big[\mu(t)^{-1}\mu^{\prime}(t)+\omega(t)\big]\sin\frac{\pi\alpha}{2}\|u^{\alpha}_t\|^{2}+\big[-2\varepsilon\omega^{\prime}(t)-4\varepsilon^{2}\omega(t)-4\varepsilon^{2}\big]\sin\frac{\pi\alpha}{2}\|u^{\alpha}\|^{2}\nonumber\\
&+4\varepsilon \sin\frac{\pi\alpha}{2}V(u^{\alpha})-8\varepsilon^{2}\mu(t)^{\frac{1-\alpha}{2}}\langle A^{\frac{1-\alpha}{4}} u^{\alpha}_t,A^{\frac{1-\alpha}{4}} u^{\alpha}\rangle.\nonumber
\end{split}
\]

Since $\frac{1+\alpha}{4}>\frac{1}{4}>\frac{1-\alpha}{4}>0$, then there are positive constants $d_0,d_1,d_2$  such that
\begin{equation}\label{1026r}
\|\cdot\|\leqslant d_0\|\cdot\|_{X^{\frac{1-\alpha}{4}}}\leqslant d_1\|\cdot\|_{_{X^{\frac{1}{4}}}}\leqslant d_2\|\cdot\|_{X^{\frac{1+\alpha}{4}}}.
\end{equation}

Let us now fix $t_0\geqslant\tau$. We restrict $\varepsilon$ in order to control $\Phi_{\alpha}$ for both sides and $\Phi_{*,\alpha}$ from below. We claim that, if we choose $\varepsilon=\varepsilon_{\omega}(t_0)=\min\Big\{1,\frac{\omega(t_0)}{4},\frac{c_1}{4(W+2)},\frac{d^{2}_{0}}{3d^{2}_{1}}\Big\}$ (as in \eqref{condE}), then 
\begin{equation}\label{Phialpha}
c_{1}\mathscr{E}_{\mathcal{H}^\alpha_t}(u^{\alpha}(t),u^{\alpha}_t(t))-2C_{\varepsilon}\Phi_{\alpha}(u^\alpha,u^\alpha_t)\leqslant \Phi_{\alpha}(u^{\alpha}(t),u^{\alpha}_t(t))\leqslant\mathscr{E}_{\mathcal{H}^\alpha_t}(u^{\alpha}(t),u^{\alpha}_t(t))+2C_{\varepsilon}
\end{equation}
and
\begin{equation}\label{1748r}
\Phi_{*,\alpha}(u^{\alpha}(t),u^{\alpha}_t(t))\geqslant  -\varepsilon\Phi_{\alpha}(u^{\alpha}(t),u^{\alpha}_t(t))
-12\varepsilon C_{\varepsilon},
\end{equation}
for all $t\in[\tau,t_0]$. Indeed, from \eqref{1559r}, we obtain
\begin{equation}\label{15r}
\begin{split}
\Phi_{\alpha}(u^\alpha,u^\alpha_t)
&\geqslant\mu(t)^{\frac{1+\alpha}{2}}\|u^\alpha\|^2_{X^{\frac{1+\alpha}{4}}}+\mu(t)^{\frac{1-\alpha}{2}}\|u_t^\alpha\|^2_{X^{\frac{1-\alpha}{4}}}
+(\omega(t)-2\varepsilon)\|u^{\alpha}\|^{2}-2C_{\varepsilon}\\
&\geqslant c_{1}\mathscr{E}_{\mathcal{H}^\alpha_t}(u^{\alpha}(t),u^{\alpha}_t(t))-2C_{\varepsilon},
\end{split}
\end{equation}
where $c_{1}=\min\{\frac{\omega(t_0)}{2},1\}$. Now, using \eqref{123}, \eqref{1026r} and Cauchy's inequality, we get 
\begin{equation}\label{16r}
\begin{split}
&\Phi_{\alpha}(u^\alpha,u^\alpha_t)\\
&\leqslant\mu(t)^{\frac{1+\alpha}{2}}\|u^\alpha\|^2_{X^{\frac{1+\alpha}{4}}}+\mu(t)^{\frac{1-\alpha}{2}}\|u_t^\alpha\|^2_{X^{\frac{1-\alpha}{4}}}
+2\varepsilon\|u^{\alpha}\|^{2}+2C_{\varepsilon}+2\varepsilon W\|u^{\alpha}\|^2\\
&+4\varepsilon\mu^{-\frac{\alpha}{2}}_{min}\mu(t)^{\frac{1+\alpha}{2}}\|u^\alpha\|^{2}_{X^{\frac{1}{4}}}+2\varepsilon \mu(t)^{\frac{1-\alpha}{2}}\|u^{\alpha}_t\|^{2}_{X^{\frac{1-\alpha}{4}}}+2\varepsilon \mu(t)^{\frac{1-\alpha}{2}}\|u^{\alpha}_t\|^{2}_{X^{\frac{1-\alpha}{4}}}\Big)\\
&\leqslant\Big(1+4\varepsilon\mu^{-\frac{\alpha}{2}}_{min}\dfrac{d^2_2}{d^2_1}+2\varepsilon\mu^{-\alpha}_{min}\dfrac{d^2_2}{d^2_0}\Big)\mu(t)^{\frac{1+\alpha}{2}}\|u^\alpha\|^2_{X^{\frac{1+\alpha}{4}}}
+2\varepsilon(1+W)\|u^{\alpha}\|^{2}\nonumber\\
&+(1+2\varepsilon)\mu(t)^{\frac{1-\alpha}{2}}\|u_t^\alpha\|^2_{X^{\frac{1-\alpha}{4}}}+2C_{\varepsilon}\\
&\leqslant c_{2}\mathscr{E}_{\mathcal{H}^\alpha_t}(u^\alpha(t),u^\alpha_t(t))+2C_{\varepsilon},
\end{split}
\end{equation}
where $c_{2}=\max\{1+4\frac{d^2_2}{d^2_1}+2\frac{d^2_2}{d^2_0},2(1+W),3\}$. Hence, using the Cauchy's inequality, \eqref{123}, \eqref{1506r}, \eqref{1559r}, \eqref{1026r}, $\omega$ is positive decreasing with \eqref{132r} and the fact that 
\[
\mathscr{E}_{\mathcal{H}^\alpha_t}(u^\alpha(t),u^\alpha_t(t))\leqslant c^{-1}_1(\Phi_{\alpha}(u^\alpha(t),u^\alpha_t(t))+2C_{\varepsilon})
\]
by using of \eqref{15r}, we obtain
\begin{equation*}
\begin{split}
&\Phi_{*,\alpha}(u^\alpha,u^\alpha_t)\\
&\geqslant\Big[-\Big(\frac{1+\alpha}{2}\Big)\mu(t)^{-1}\mu'(t)-6\varepsilon\Big]\mu(t)^{\frac{1-\alpha}{2}} \|u^{\alpha}_t\|^{2}_{X^{\frac{1-\alpha}{4}}}+4\cos\frac{\pi\alpha}{2}\mu(t)^{\frac{1}{2}}\|u^{\alpha}_t\|^{2}_{X^{\frac{1}{4}}}\\
&+\big[\mu(t)^{-1}\mu^{\prime}(t)+\omega(t)\big]\sin\frac{\pi\alpha}{2}\|u^{\alpha}_t\|^{2}-4\varepsilon C_{\varepsilon}\\
&+\Big[-\Big(\frac{1+\alpha}{2}\Big)\mu(t)^{-1}\mu'(t)+2\varepsilon\Big]\mu(t)^{\frac{1+\alpha}{2}} \|u^{\alpha}\|^{2}_{X^{\frac{1+\alpha}{4}}}\\
&-8\varepsilon^{2}\cos\frac{\pi\alpha}{2}\mu(t)^{\frac{1}{2}}\|u^{\alpha}\|^{2}_{X^{\frac{1}{4}}}+\big[-2\varepsilon\omega^{\prime}(t)-4\varepsilon^{2}\omega(t)-8\varepsilon^{2}\big]\sin\frac{\pi\alpha}{2}\|u^{\alpha}\|^{2}\\
&-4\varepsilon^{2}\mu(t)^{\frac{1-\alpha}{2}}\|u^{\alpha}_t\|^{2}_{X^{\frac{1-\alpha}{4}}} -4\varepsilon^{2}\mu(t)^{\frac{1-\alpha}{2}}\|u^{\alpha}\|^{2}_{X^{\frac{1-\alpha}{4}}}\\
&\geqslant(-\mu(t)^{-1}\mu'(t)-4\varepsilon^{2})\mu(t)^{\frac{1-\alpha}{2}} \|u^{\alpha}_t\|^{2}_{X^{\frac{1-\alpha}{4}}}-4\varepsilon C_{\varepsilon}\\
&+4\cos\frac{\pi\alpha}{2}\mu(t)^{\frac{1}{2}}\|u^{\alpha}_t\|^{2}_{X^{\frac{1}{4}}}-6\varepsilon\mu(t)^{\frac{1-\alpha}{2}} \|u^{\alpha}_t\|^{2}_{X^{\frac{1-\alpha}{4}}}\\
&+(-\mu(t)^{-1}\mu'(t)+2\varepsilon)\mu(t)^{\frac{1+\alpha}{2}} \|u^{\alpha}\|^{2}_{X^{\frac{1+\alpha}{4}}}\\
&-8\varepsilon^{2}\cos\frac{\pi\alpha}{2}\mu(t)^{\frac{1}{2}}\|u^{\alpha}\|^{2}_{X^{\frac{1}{4}}}-4\varepsilon^{2}(W+2)\|u^{\alpha}\|^{2}-4\varepsilon^{2}\mu(t)^{\frac{1-\alpha}{2}}\|u^{\alpha}\|^{2}_{X^{\frac{1-\alpha}{4}}}\\
&\geqslant-5\varepsilon^{2}\mu(t)^{\frac{1-\alpha}{2}} \|u^{\alpha}_t\|^{2}_{X^{\frac{1-\alpha}{4}}}+\Big(4\mu^{\alpha}_{min}\cos\frac{\pi\alpha}{2}-6\dfrac{d^{2}_{1}}{d^{2}_{0}}\varepsilon\Big)\mu(t)^{\frac{1-\alpha}{2}} \|u^{\alpha}_t\|^{2}_{X^{\frac{1}{4}}}\\
&-4\varepsilon^{2}\Big(2\cos\frac{\pi\alpha}{2}\mu^{-\frac{\alpha}{2}}_{min}\dfrac{d^{2}_{2}}{d^{2}_{1}}+\mu^{-\alpha}_{min}\dfrac{d^{2}_{2}}{d^{2}_{0}}\Big)\mu(t)^{\frac{1+\alpha}{2}}\|u^{\alpha}\|^{2}_{X^{\frac{1+\alpha}{4}}}-4\varepsilon^{2}(W+2)\|u^{\alpha}\|^{2}-4\varepsilon C_{\varepsilon}\\
&\geqslant  -4\varepsilon^{2}(W+2)\mathscr{E}_{\mathcal{H}^\alpha_t}(u^\alpha(t),u^\alpha_t(t))-4\varepsilon C_{\varepsilon}\\
&\geqslant  -\varepsilon\Phi_{\alpha}(t)
-12\varepsilon C_{\varepsilon},
\end{split}
\end{equation*}
as claimed. Here, we have used that there exists $\alpha_0<1$ such that for all $\alpha\in(\alpha_0,1)$
\[
2\cos\frac{\pi\alpha}{2}\mu^{-\frac{\alpha}{2}}_{min}\dfrac{d^{2}_{2}}{d^{2}_{1}}\leqslant \dfrac{W+2}{2}
\]
and $W$ is chosen large enough so that 
\[
\mu^{-\alpha}_{min}\dfrac{d^{2}_{2}}{d^{2}_{0}}\leqslant \dfrac{W+2}{2}.
\]

Thus, using \eqref{10r} and \eqref{1748r}, we get 
\begin{equation}\label{11r}
\dfrac{d}{dt}\Phi_{\alpha}(u^\alpha(t),u^\alpha_t(t))+\varepsilon\Phi_{\alpha}(u^\alpha(t),u^\alpha_t(t))\leqslant 16\varepsilon C_{\varepsilon}.
\end{equation}
Multiplying \eqref{11r} by $e^{\varepsilon t}$ and integrating between $\tau$ and $t_0$, we obtain
\begin{equation}\label{12r}
\Phi_{\alpha}(u^\alpha(t_0),u^\alpha_t(t_0))\leqslant\Phi_{\alpha}(u^\alpha(\tau),u^\alpha_t(\tau))e^{-\varepsilon(t_0-\tau)} +16C_{\varepsilon}.
\end{equation}
Now, from \eqref{15r} and \eqref{16r}, we obtain
\begin{equation}\label{17r}
\begin{split}
c_1\mathscr{E}_{\mathcal{H}^\alpha_t}(u^\alpha(t_0),u^\alpha_t(t_0))
&\leqslant\Phi_{\alpha}(u^\alpha(t_0),u^\alpha_t(t_0))+2C_{\varepsilon}\\
&\leqslant \Phi_{\alpha}(u^\alpha(\tau),u^\alpha_t(\tau))e^{-\varepsilon(t_0-\tau)}+18C_{\varepsilon}\\
&\leqslant c_2\mathscr{E}_{\mathcal{H}^\alpha_\tau}(u^\alpha(\tau),u^\alpha_t(\tau))e^{-\varepsilon(t_0-\tau)}+20C_{\varepsilon}.
\end{split}
\end{equation}
Furthermore,
\[
\mathscr{E}_{\mathcal{H}^\alpha_t}(u^\alpha(t),u^\alpha_t(t))\leqslant M_0\mathscr{E}_{\mathcal{H}^\alpha_s}(u^\alpha(s),u^\alpha_t(s))e^{-\varepsilon(t-\tau)}+M_1,
\]
where $M_0=c_2/c_1$ and $M_1=20C_{\varepsilon}/c_1$.   \qed

\bigskip

 We recall that the function $\Phi_{\alpha}(u^{\alpha},v^{\alpha})$ (see \eqref{funcional01}) along with 
\[
v^\alpha=\mu(t)^{\frac{-1+\alpha}{2}}\sin\frac{\pi\alpha}{2}A^{\frac{-1+\alpha}{2}}v^\alpha-\mu(t)^{\frac{\alpha}{2}} \cos\frac{\pi\alpha}{2}A^{\frac{\alpha}{2}}u^\alpha
\]
 (see \eqref{6f}), define a functional $\mathscr{L}_{\alpha}$ by
\begin{equation}\label{funcionalL}
\begin{split}
& \mathscr{L}_{\alpha}(\left[\begin{smallmatrix}u^{\alpha}\\ v^{\alpha} \end{smallmatrix}\right]) \\ 
&=\mu(t)^{\frac{1+\alpha}{2}}\|u^\alpha\|^2_{X^{\frac{1+\alpha}{4}}}+4\varepsilon\cos\frac{\pi\alpha}{2}\mu(t)^{\frac{1}{2}} \|u^\alpha\|^2_{X^{\frac{1}{4}} }+2\varepsilon\sin\frac{\pi\alpha}{2}\omega(t)\|u^\alpha\|^2\\
&+\mu(t)^{\frac{1-\alpha}{2}}\|\mu(t)^{\frac{-1+\alpha}{2}}\sin\frac{\pi\alpha}{2}A^{\frac{-1+\alpha}{2}}v^\alpha-\mu(t)^{\frac{\alpha}{2}} \cos\frac{\pi\alpha}{2}A^{\frac{\alpha}{2}}u^\alpha\|^2_{X^{\frac{1-\alpha}{4}}}-2\sin\frac{\pi\alpha}{2}V(u^{\alpha})\\
&+4\varepsilon\mu(t)^{\frac{1-\alpha}{2}}\langle A^{\frac{1-\alpha}{4}}(\mu(t)^{\frac{-1+\alpha}{2}}\sin\frac{\pi\alpha}{2}A^{\frac{-1+\alpha}{2}}v^\alpha-\mu(t)^{\frac{\alpha}{2}} \cos\frac{\pi\alpha}{2}A^{\frac{\alpha}{2}}u^\alpha), A^{\frac{1-\alpha}{4}} u^{\alpha}\rangle\\
&=\mu(t)^{\frac{1+\alpha}{2}}\|u^\alpha\|^2_{X^{\frac{1+\alpha}{4}}}+\|\sin\frac{\pi\alpha}{2}A^{\frac{-1+\alpha}{4}}v^\alpha-\mu(t)^{\frac{1}{2}} \cos\frac{\pi\alpha}{2}A^{\frac{1+\alpha}{4}}u^\alpha\|^2-2\sin\frac{\pi\alpha}{2}V(u^{\alpha})\\
&+2\varepsilon\sin\frac{\pi\alpha}{2}\omega(t)\|u^\alpha\|^2+4\varepsilon\sin\frac{\pi\alpha}{2}\langle A^{\frac{1-\alpha}{4}}u^\alpha, A^{\frac{-1+\alpha}{4}} v^{\alpha}\rangle,
\end{split}
\end{equation}
with domain  
\begin{equation}\label{funcionalLD}
D(\mathscr{L}_{\alpha})=\Big\{\left[\begin{smallmatrix}u^{\alpha}\\ v^{\alpha} \end{smallmatrix}\right]\in X^\frac{1+\alpha}{4}\times X^\frac{-1+\alpha}{4}: \int_0^{u^{\alpha}}f(\sigma)d\sigma\in L^{1}(\Omega)\Big\}. 
\end{equation}

\begin{remark}
On the functional $\mathscr{L}_{\alpha}$ defined in \eqref{funcionalL}-\eqref{funcionalLD}, we have
\begin{itemize}
\item[$i)$] Observe that $D(\mathscr{L}_{\alpha})= X^\frac{1+\alpha}{4}\times X^\frac{-1+\alpha}{4}$ provided that $\alpha$ is close enough to $1$.
\item[$ii)$] In particular, if $s\in(\frac{N}2(1-\frac1\rho) -\frac1\rho,1)$ then, due to \eqref{eq:regularity-of-solution}, $\mathscr{L}_{\alpha}(\left[\begin{smallmatrix}u^{\alpha}\\ v^{\alpha} \end{smallmatrix}\right])$ is well defined for all $\alpha\in[\frac{N}2(\rho-1) -\rho s,1)$ and $t\in(\tau,\xi_{u_{\tau\alpha},v_{\tau\alpha}})$ along each solution $\left[\begin{smallmatrix}u^{\alpha}\\ v^{\alpha} \end{smallmatrix}\right]$ through $\left[\begin{smallmatrix}u_{\tau\alpha}\\ v_{\tau\alpha} \end{smallmatrix}\right]\in X^\frac{s}{2}\times X^\frac{s-1}{2}$ from Theorem \ref{th:local-well-posedness-of-perturbed-problem}.
\item[$iii)$] In fact, for positive times and as long as the solutions exist we have
\begin{equation}\label{2216r}
\Phi_{\alpha}(u^{\alpha},v^{\alpha})=\mathscr{L}_{\alpha}\Big(\left[\begin{smallmatrix}u^{\alpha}(t)\\ v^{\alpha}(t) \end{smallmatrix}\right]\Big) \quad\text{ and }\quad\dfrac{d}{dt}(\Phi_{\alpha}(u^{\alpha},v^{\alpha}))=\dfrac{d}{dt}\Big(\mathscr{L}_{\alpha}\Big(\left[\begin{smallmatrix}u^{\alpha}(t)\\ v^{\alpha}(t) \end{smallmatrix}\right]\Big) \Big).
\end{equation}
\end{itemize}
\end{remark}

The next result shows the equivalence between the energy $\mathscr{E}_{\mathcal{H}^\alpha_t}(\cdot,\cdot)$ and the norm of the space $X^{\frac{1+\alpha}{4}}\times X^{\frac{-1+\alpha}{4}}$ due to condition \eqref{123} and the form of the system \eqref{6f}.

\begin{lemma}\label{equiv-geral}
Under the same conditions as in Lemma \ref{inequality-geral}, we have
\begin{equation}\label{Equivar}
C_1\|\left[\begin{smallmatrix}u^{\alpha}\\ v^{\alpha}\end{smallmatrix}\right]\|_{X^{\frac{1+\alpha}{4}}\times X^{\frac{-1+\alpha}{4}}}\leqslant\mathscr{E}_{\mathcal{H}^\alpha_t}(u^{\alpha},v^{\alpha})\leqslant C_2 \|\left[\begin{smallmatrix}u^{\alpha}\\ v^{\alpha}\end{smallmatrix}\right]\|_{X^{\frac{1+\alpha}{4}}\times X^{\frac{-1+\alpha}{4}}}, \quad\forall t\geqslant \tau,
\end{equation}
 where $C_1$ and $C_2$ are positive constants.
\end{lemma}

\proof Firstly, we see that $\frac{1+\alpha}{4}>\frac{1-\alpha}{4}>0>\frac{-1+\alpha}{4}$, then there are positive constants $d_3,d_4,d_5$  such that
\begin{equation}\label{1026ro}
\|\cdot\|_{X^{\frac{-1+\alpha}{4}}}\leqslant d_{3}\|\cdot\|\leqslant d_4\|\cdot\|_{X^{\frac{1-\alpha}{4}}}\leqslant d_5\|\cdot\|_{X^{\frac{1+\alpha}{4}}}.
\end{equation}
Hence, using the first equation of system \eqref{6f} and \eqref{1026ro}, we obtain
\[
\begin{split}
\mathscr{E}_{\mathcal{H}^\alpha_t}(u^{\alpha},v^{\alpha})&\leqslant \Big(\mu_{max}+\dfrac{d^2_5}{d^2_3}\Big)\|u^{\alpha}\|^2_{X^{\frac{1+\alpha}{4}}}\\
&+\mu(t)^{\frac{1-\alpha}{2}}\|\mu(t)^{\frac{-1+\alpha}{2}}\sin\frac{\pi\alpha}{2}A^{\frac{-1+\alpha}{2}}v^{\alpha}-\mu(t)^{\frac{\alpha}{2}} \cos\frac{\pi\alpha}{2}A^{\frac{\alpha}{2}}u^{\alpha}\|^2_{X^{\frac{1-\alpha}{4}}}\\
&\leqslant \Big(\mu_{max}+\dfrac{d^2_5}{d^2_3}\Big)\|u^{\alpha}\|^2_{X^{\frac{1+\alpha}{4}}}+2\sin^2\frac{\pi\alpha}{2}\|v^{\alpha}\|^2_{X^{\frac{-1+\alpha}{4}}}+2\mu(t) \cos^2\frac{\pi\alpha}{2}\|u^{\alpha}\|^2_{X^{\frac{1+\alpha}{4}}}\nonumber\\
&\leqslant C_{2}(\|u^{\alpha}\|^2_{X^{\frac{1+\alpha}{4}}}+\|v^{\alpha}\|^2_{X^{\frac{-1+\alpha}{4}}}),
\end{split}
\]
where $C_2=\max\{3\mu_{max}+\frac{d^2_5}{d^2_3},2\}$. Finally, we have
\[
\begin{split}
\mathscr{E}_{\mathcal{H}^\alpha_t}(u^{\alpha},v^{\alpha})&\geqslant \mu(t)^{\frac{1+\alpha}{2}}\|u^\alpha\|^2_{X^{\frac{1+\alpha}{4}}}+\mu(t)^{\frac{1-\alpha}{2}}\|v^\alpha\|^2_{X^{\frac{1-\alpha}{4}}}\geqslant C_{1}(\|u^{\alpha}\|^2_{X^{\frac{1+\alpha}{4}}}+\|v^{\alpha}\|^2_{X^{\frac{-1+\alpha}{4}}}),
\end{split}
\]
where $C_1=\min\{1,\frac{1}{d^2_4}\}$.
\qed

\medskip

We now prove that the functional $\mathscr{L}_{\alpha}(\cdot)$ is bounded from both sides as stated in the following lemma.

\begin{lemma}
Suppose the same conditions as in Lemma \ref{inequality-geral}. Then, there are positive constants $D_1,D_2,D_3$ such that for all $\alpha <1$ close enough to $1$, $\mathscr{L}_{\alpha}$ in \eqref{funcionalL} satisfies
\begin{equation}\label{equivaLR}
D_1\Big\|\left[\begin{smallmatrix}u^{\alpha}(t)\\ v^{\alpha}(t)\end{smallmatrix}\right]\Big\|_{X^{\frac{1+\alpha}{4}}\times X^{\frac{-1+\alpha}{4}}}-D_3\leqslant\mathscr{L}_{\alpha}\Big(\left[\begin{smallmatrix}u^{\alpha}(t)\\ v^{\alpha}(t) \end{smallmatrix}\right]\Big)\leqslant D_2 \Big\|\left[\begin{smallmatrix}u^{\alpha}(t)\\ v^{\alpha}(t)\end{smallmatrix}\right]\Big\|_{X^{\frac{1+\alpha}{4}}\times X^{\frac{-1+\alpha}{4}}}+D_3, 
\end{equation}
for all $\left[\begin{smallmatrix}u^{\alpha}\\ v^{\alpha}\end{smallmatrix}\right]\in X^{\frac{1+\alpha}{4}}\times X^{\frac{-1+\alpha}{4}}$, and all $t\in\mathbb{R}$.
\end{lemma}

\proof The result follows directly from Lemma \ref{equiv-geral}, and \eqref{Phialpha} where $D_1=c_1 C_1$, $D_2=c_2 C_2$ and $D_3=2 C_{\varepsilon}$.
\qed

\subsection{Existence of absorbing and attractor}

In this subsection, we show the main results of this paper, the existence of a pullback absorbing family to bounded sets and existence of pullback attractor for the  Cauchy problem \eqref{ffhghgf}.

\begin{theorem}\label{theo:10r}
Suppose that $\mu$ is $(\gamma, \kappa)-$H\"{o}lder continuous in $\mathbb{R}$ satisfying \eqref{123} and $E^{\theta}_t(\alpha)$ as in \eqref{eq:characterization}. If $s\in(\frac{N}2(1-\frac1\rho) -\frac1\rho,1)$, then, for any
$\alpha\in[\frac{N}2(\rho-1) - \rho s,1)$, there exist $\theta_1 = \frac{-\frac{N}{2}(\rho - 1) + \rho s}{\alpha}\in[-1, 0)$ and $\theta_2 = \frac{s-1}{\alpha}\in(\theta_1,1)$ such that $0<\theta_2 - \theta_1 < 1$ and the following hold:
\begin{itemize}
\item[$i)$] The operators $\varLambda(t)^{\alpha}$ are uniformly sectorial and the map $t \mapsto \varLambda(t)^{\alpha}$ is uniformly H\"{o}lder continuous in $E^{\theta_1}_t(\alpha)$;
\item[$ii)$] There exists a linear evolution process $\{U_{\alpha}(t,\tau): t\geqslant\tau\in\mathbb{R}\}$ that solves the linear homogeneous problem 
\begin{equation}\label{128r}
\begin{cases}
\dfrac{d}{dt}\left[\begin{smallmatrix} u^{\alpha}\\ v^{\alpha} \end{smallmatrix}\right]+\varLambda(t)^\alpha \left[\begin{smallmatrix} u^{\alpha}\\ v^{\alpha} \end{smallmatrix}\right]=0,\quad t>\tau,\\ \vspace{-0.4cm}\\
\left[\begin{smallmatrix} u^{\alpha}\\ v^{\alpha} \end{smallmatrix}\right](\tau)=\left[\begin{smallmatrix} u_{\tau\alpha}\\ v_{\tau\alpha} \end{smallmatrix}\right]\in E^{\theta_2}_\tau(\alpha),
\end{cases}
\end{equation}
that is, for any $t\geqslant\tau\in\mathbb{R}$,
\[
U_{\alpha}(t,\tau)\left[\begin{smallmatrix} u_{\tau\alpha}\\ v_{\tau\alpha} \end{smallmatrix}\right]=\left[\begin{smallmatrix} u^{\alpha}\\ v^{\alpha} \end{smallmatrix}\right],
\]
is given by
\[
U_\alpha(t,\tau) = e^{-(t-\tau)\varLambda(\tau)^\alpha} + \int_{\tau}^{t}U_{\alpha}(t,s)[\varLambda(\tau)^\alpha -\varLambda(s)^\alpha]e^{-(s-\tau)\varLambda(\tau)^\alpha}ds,
\]
where
$$(\tau,\infty)\ni t\mapsto \left[\begin{smallmatrix} u^{\alpha}\\ v^{\alpha} \end{smallmatrix}\right](t)= U_{\alpha}(t,\tau)\left[\begin{smallmatrix} u_{\tau\alpha}\\ v_{\tau\alpha} \end{smallmatrix}\right] \in E^{\theta_2}_t(\alpha)$$ is continuously differentiable,
$$\left[\begin{smallmatrix} u^{\alpha}\\ v^{\alpha}\end{smallmatrix}\right](t)\in X^{\frac{1+\alpha}{2}}\times X^{\frac{\alpha}{2}},\quad \forall t\in (\tau,\infty)$$
 and satisfies \eqref{128r}.
 \end{itemize}
 \end{theorem}

\proof 
First we claimed that $\varLambda(t)^{\alpha}$ is uniformly sectorial in $E^{\theta_1}_t(\alpha)$; that is, there is a constant $C > 0$ such that
\begin{equation}\label{123r}
\|(\lambda I+\varLambda(t)^{\alpha})^{-1}\|_{\mathscr{L}(E^{\theta_1}_t(\alpha))}\leqslant\dfrac{C}{|\lambda|+1}, \quad\text{ for all }\lambda\in\mathbb{C}\text{ with  Re}\,\lambda\geqslant 0.
\end{equation} 
This fact follows from sectoriallity of the operators $A(t)$, \eqref{123}, Amann \cite[Theorem 2.1.3]{Am} and using the ideas of  Bezerra and Nascimento \cite[Theorem 1.1]{BN}.

 By the other hand, by using item $vi)$ of Lemma \ref{lem:useful:properties}, we can see that for any $t,\tau,s\in\mathbb{R}$,
\[
[\varLambda(t)^{\alpha}-\varLambda(\tau)^{\alpha}]\varLambda(s)^{-\alpha}=
\begin{bmatrix} 
E_{11} & E_{12} \\
E_{21} & E_{22} 
\end{bmatrix},
\]
where
\[
\begin{split}
E_{11}&= \cos^{2}\frac{\pi\alpha}{2}[\mu(t)^{\frac{\alpha}{2}} - \mu(\tau)^{\frac{\alpha}{2}}]\mu(s)^{-\frac{\alpha}{2}}I + \sin^{2}\frac{\pi\alpha}{2}[\mu(t)^{\frac{-1+\alpha}{2}} - \mu(\tau)^{\frac{-1+\alpha}{2}}]\mu(s)^{\frac{1-\alpha}{2}}I,\\
E_{12}&=\cos\frac{\pi\alpha}{2}\sin\frac{\pi\alpha}{2}\Big\{[\mu(t)^{\frac{\alpha}{2}} - \mu(\tau)^{\frac{\alpha}{2}}]\mu(s)^{\frac{-1-\alpha}{2}} - [\mu(t)^{\frac{-1+\alpha}{2}} - \mu(\tau)^{\frac{-1+\alpha}{2}}]\mu(s)^{-\frac{\alpha}{2}}\Big\}A^{-\frac{1}{2}},\\
E_{21}&=\cos\frac{\pi\alpha}{2}\sin\frac{\pi\alpha}{2}\Big\{[\mu(t)^{\frac{1+\alpha}{2}} - \mu(\tau)^{\frac{1+\alpha}{2}}]\mu(s)^{-\frac{\alpha}{2}} - [\mu(t)^{\frac{\alpha}{2}} - \mu(\tau)^{\frac{\alpha}{2}}]\mu(s)^{\frac{1-\alpha}{2}}\Big\}A^{\frac{1}{2}},\\
E_{11}&= \sin^{2}\frac{\pi\alpha}{2}[\mu(t)^{\frac{1+\alpha}{2}} - \mu(\tau)^{\frac{1+\alpha}{2}}]\mu(s)^{\frac{-1-\alpha}{2}}I+ \cos^{2}\frac{\pi\alpha}{2}[\mu(t)^{\frac{\alpha}{2}} - \mu(\tau)^{\frac{\alpha}{2}}]\mu(s)^{-\frac{\alpha}{2}}I.
\end{split}
\]
Thus, 
\[
\begin{split}
\Big\|&[\varLambda(t)^{\alpha}-\varLambda(\tau)^{\alpha}]\varLambda(s)^{-\alpha}
\begin{bmatrix} 
u^{\alpha} \\
v^{\alpha} 
\end{bmatrix}
\Big\|_{E^{\theta_1}_t(\alpha)}
=
\Big\|
\begin{bmatrix} 
E_{11}u^{\alpha} + E_{12}v^{\alpha} \\
E_{21}u^{\alpha} + E_{22}v^{\alpha} 
\end{bmatrix}
\Big\|_{X^{\frac{1+\alpha\theta_1}{2}}\times X^{\frac{\alpha\theta_1}{2}}}\\
&=\|A^{\frac{1+\alpha\theta_1}{2}}[E_{11}u^{\alpha} + E_{12}v^{\alpha}]\|+\|A^{\frac{\alpha\theta_1}{2}}[E_{21}u^{\alpha} + E_{22}v^{\alpha}]\|\\
&\leqslant \max\{1,\mu^{\frac{1}{2}}_{max}\}\mu^{-\frac{\alpha}{2}}_{min}
[|\mu(t)^{\frac{\alpha}{2}} - \mu(\tau)^{\frac{\alpha}{2}}| + |\mu(t)^{\frac{-1+\alpha}{2}} - \mu(\tau)^{\frac{-1+\alpha}{2}}|]
\|A^{\frac{1+\alpha\theta_1}{2}}u^{\alpha}\|\\
&+\max\{1,\mu^{-\frac{1}{2}}_{max}\}\mu^{-\frac{\alpha}{2}}_{min}
[|\mu(t)^{\frac{\alpha}{2}} - \mu(\tau)^{\frac{\alpha}{2}}| + |\mu(t)^{\frac{-1+\alpha}{2}} - \mu(\tau)^{\frac{-1+\alpha}{2}}|]
\|A^{\frac{\alpha\theta_1}{2}}v^{\alpha}\|\\
&+\mu^{-\frac{\alpha}{2}}_{min}[|\mu(t)^{\frac{1+\alpha}{2}} - \mu(\tau)^{\frac{1+\alpha}{2}}|+|\mu(t)^{\frac{\alpha}{2}} - \mu(\tau)^{\frac{\alpha}{2}}|]\|A^{\frac{1+\alpha\theta_1}{2}}u^{\alpha}\|\\
&+\max\{1,\mu^{-\frac{1}{2}}_{min}\}\mu^{-\frac{\alpha}{2}}_{min}
[|\mu(t)^{\frac{1+\alpha}{2}} - \mu(\tau)^{\frac{1+\alpha}{2}}| + |\mu(t)^{\frac{\alpha}{2}} - \mu(\tau)^{\frac{\alpha}{2}}|]
\|A^{\frac{\alpha\theta_1}{2}}v^{\alpha}\|\\
&\leqslant 2\max\{\mu^{-\frac{1}{2}}_{max},\mu^{\frac{1}{2}}_{max}\}\mu^{-\frac{\alpha}{2}}_{min}
[|\mu(t)^{\frac{\alpha}{2}} - \mu(\tau)^{\frac{\alpha}{2}}| + |\mu(t)^{\frac{-1+\alpha}{2}} - \mu(\tau)^{\frac{-1+\alpha}{2}}|]
\Big\|
\begin{bmatrix} 
u^{\alpha} \\
v^{\alpha}
\end{bmatrix}
\Big\|_{E^{\theta_1}_t(\alpha)}\\
&+2\max\{1,\mu^{-\frac{1}{2}}_{min}\}\mu^{-\frac{\alpha}{2}}_{min}
[|\mu(t)^{\frac{1+\alpha}{2}} - \mu(\tau)^{\frac{1+\alpha}{2}}| + |\mu(t)^{\frac{\alpha}{2}} - \mu(\tau)^{\frac{\alpha}{2}}|]
\Big\|
\begin{bmatrix} 
u^{\alpha} \\
v^{\alpha}
\end{bmatrix}
\Big\|_{E^{\theta_1}_t(\alpha)}.
\end{split}
\]
Since, for all $t,\tau\in\mathbb{R}$,  
\[
\begin{split}
|\mu(t)^{\frac{-1+\alpha}{2}} - \mu(\tau)^{\frac{-1+\alpha}{2}}| &\leqslant \mu^{-\frac{1}{2}}_{min}|\mu(t)^{\frac{\alpha}{2}} - \mu(\tau)^{\frac{\alpha}{2}}|+\max\{1,\mu^{\frac{1}{2}}_{max}\}|\mu(t)^{\frac{1}{2}} - \mu(\tau)^{\frac{1}{2}}|,\\
|\mu(t)^{\frac{1+\alpha}{2}} - \mu(\tau)^{\frac{1+\alpha}{2}}| &\leqslant \mu^{\frac{1}{2}}_{max}|\mu(t)^{\frac{\alpha}{2}} - \mu(\tau)^{\frac{\alpha}{2}}|+\max\{1,\mu^{\frac{1}{2}}_{max}\}|\mu(t)^{\frac{1}{2}} - \mu(\tau)^{\frac{1}{2}}|.
\end{split}
\]
Then, there is a positive constant $C^{\prime}=C^{\prime}(\mu_{min},\mu_{max})$ (independent of $\alpha$) such that 
\[
\begin{split}
\Big\|&[\varLambda(t)^{\alpha}-\varLambda(\tau)^{\alpha}]\varLambda(s)^{-\alpha}
\begin{bmatrix} 
u^{\alpha} \\
v^{\alpha} 
\end{bmatrix}
\Big\|_{E^{\theta_1}_t(\alpha)}\leqslant C^{\prime}[|\mu(t)^{\frac{\alpha}{2}} - \mu(\tau)^{\frac{\alpha}{2}}| + |\mu(t)^{\frac{1}{2}} - \mu(\tau)^{\frac{1}{2}}|]
\Big\|
\begin{bmatrix} 
u^{\alpha} \\
v^{\alpha}
\end{bmatrix}
\Big\|_{E^{\theta_1}_t(\alpha)}.
\end{split}
\]

From Bezerra and Nascimento \cite[Lemma 3.6]{BN}, $\mu^{\frac{\alpha}{2}}$ is $(\frac{1}{4},\kappa)-$H\"{o}lder conti\-nuous in $\mathbb{R}$, so that there exists a constant $C= C(\mu_{min},\mu_{max},\kappa)>0$  (independent of $\alpha$) such that 
\begin{equation}\label{126r}
\Big\|[\varLambda(t)^{\alpha}-\varLambda(\tau)^{\alpha}]\varLambda(s)^{-\alpha}
\begin{bmatrix} 
u^{\alpha} \\
v^{\alpha} 
\end{bmatrix}
\Big\|_{E^{\theta_1}_t(\alpha)}\leqslant C|t-\tau|^{\frac{1}{4}}
\Big\|
\begin{bmatrix} 
u^{\alpha} \\
v^{\alpha}
\end{bmatrix}
\Big\|_{E^{\theta_1}_t(\alpha)}, \forall t,\tau\in\mathbb{R}. 
\end{equation}

Item $i)$ is proved using \eqref{123r} and \eqref{126r}, and item $ii)$ follows from Carvalho and Nascimento \cite[Section 2]{CN} and Sobolevski\u{\i} \cite{sobol}, .
\qed

\medskip

For a better understanding of the relationship of the scale of fractional power spaces of operator $\varLambda(t)$, we construct the following diagram:
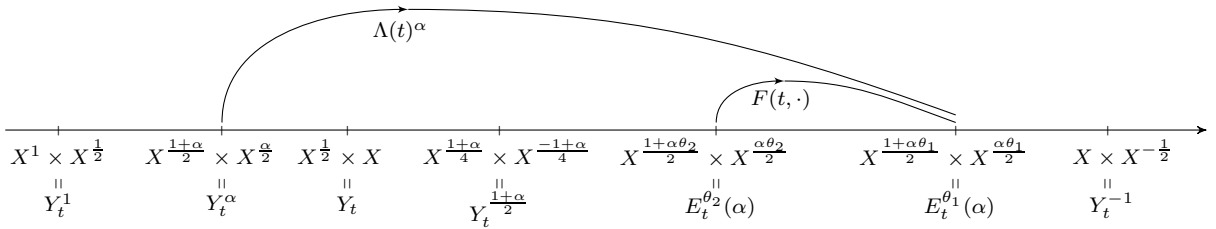
\begin{figure}[h]
\begin{center}
\begin{tikzpicture}
\draw[-stealth'] (-4.5,0) -- (11.3,0);
\node at (-3.8,0){\tiny $\shortmid$};
\node at (-1.65,0){\tiny $\shortmid$};
\node at (0,0){\tiny $\shortmid$};
\node at (2,0){\tiny $\shortmid$};
\node at (4.85,0){\tiny $\shortmid$};
\node at (8,0){\tiny $\shortmid$};
\node at (10,0){\tiny $\shortmid$};

\node at (0.7,1.3) {\tiny $\Lambda(t)^{\alpha}$};
\node at (5.7,0.4) {\tiny $F(t,\cdot)$};

\node at (-3.8,-0.3){\tiny $X^{1}\times X^{\frac{1}{2}}$};
\node at (-3.8,-0.7){\tiny $\shortparallel$};
\node at (-3.8,-1){\tiny $Y^{1}_t$};

\node at (-1.8 ,-0.3){\tiny $X^{\frac{1+\alpha}{2}}\times X^{\frac{\alpha}{2}}$};
\node at (-1.65,-0.7){\tiny $\shortparallel$};
\node at (-1.65,-1){\tiny $Y^{\alpha}_t$};

\node at (-0.1,-0.3){\tiny $X^{\frac{1}{2}}\times X$};
\node at (0,-0.7){\tiny $\shortparallel$};
\node at (0,-1){\tiny $Y_t$} ;

\node at (2.1,-0.3){\tiny $X^{\frac{1+\alpha}{4}}\times X^{\frac{-1+\alpha}{4}}$};
\node at (2,-0.7){\tiny $\shortparallel$};
\node at (2,-1.05) {\tiny $Y^{\frac{1+\alpha}{2}}_t$};

\node at (4.7,-0.3){\tiny $X^{\frac{1+\alpha\theta_2}{2}}\times X^{\frac{\alpha\theta_2}{2}}$};
\node at (4.85,-0.7){\tiny $\shortparallel$};
\node at (4.9,-1){\tiny $E^{\theta_2}_t(\alpha)$} ;

\node at (7.85,-0.3){\tiny $X^{\frac{1+\alpha\theta_1}{2}}\times X^{\frac{\alpha\theta_1}{2}}$};
\node at (8,-0.7){\tiny $\shortparallel$};
\node at (8.05,-1) {\tiny $E^{\theta_1}_t(\alpha)$};

\node at (10.2,-0.3){\tiny $X\times X^{-\frac{1}{2}}$};
\node at (10,-0.7){\tiny $\shortparallel$};
\node at (10.05,-1){\tiny $Y^{-1}_t$};
\draw[-latex'] (-1.65,0.1) to [out=90,in=180] (0.8,1.6);
\draw (0.8,1.6) to [out=0,in=160] (8,0.2);
\draw[-latex'] (4.85,0.1) to [out=90,in=180] (5.75,0.65);
\draw (5.75,0.65) to [out=0,in=160] (8,0.1);
\end{tikzpicture}
\caption{\small Partial description of the fractional power spaces scale for $\varLambda(t)$, $t\in\mathbb{R}$.}\label{figuraScalar}
\end{center}
\end{figure}

\begin{lemma}\label{lemma:CD02}
Suppose that \eqref{growth-condition} and \eqref{eq:supercritical} hold, $\omega$ is $(\zeta, \kappa_0)-$H\"{o}lder continuous in $\mathbb{R}$ satisfying \eqref{132r} and $E^{\theta}_t(\alpha)$ as in \eqref{eq:characterization}. If $\frac{N}{2}(1-\frac{1}{\rho})-\frac{1}{\rho} < s \leq\frac{N}{2}(1-\frac{1}{\rho})$, then for any $\alpha\in[s - 1 + \frac{N}{2} (\rho - 1) - \rho s, 1)$ satisfying $\frac{N}{2}(\rho - 1) - \rho s < \alpha$, there exist $\theta_1 = \frac{-\frac{N}{2}(\rho - 1) + \rho s}{\alpha}\in[-1, 0)$ and $\theta_2 = \frac{s-1}{\alpha}\in(\theta_1,1)$ such that $0<\theta_2 - \theta_1 < 1$ and for $F$ in \eqref{eq:Lambdaf}, we have
\begin{equation}\label{140r}
\|F(t,\left[\begin{smallmatrix} u^{\alpha}\\ v^{\alpha} \end{smallmatrix}\right])\|_{E^{\theta_1}_t(\alpha)}\leqslant c(1+\|\left[\begin{smallmatrix} u^{\alpha}\\ v^{\alpha} \end{smallmatrix}\right]\|^{\rho}_{E^{\theta_2}_t(\alpha)}), \quad \left[\begin{smallmatrix} u^{\alpha}\\ v^{\alpha} \end{smallmatrix}\right] \in E^{\theta_2}_t(\alpha)
\end{equation}
and 
\[
\begin{split}
\Big\|F\Big(t,\left[\begin{smallmatrix} u^{\alpha}_1\\ v^{\alpha}_1 \end{smallmatrix}\right]\Big)-F\Big(\tau,\left[\begin{smallmatrix} u^{\alpha}_2\\ v^{\alpha}_2 \end{smallmatrix}\right]\Big)\Big\|_{E^{\theta_1}_t(\alpha)}\leqslant &c\Big(|t-\tau|^{\zeta}
+\Big\|\left[\begin{smallmatrix} u^{\alpha}_1\\ v^{\alpha}_1 \end{smallmatrix}\right]-\left[\begin{smallmatrix} u^{\alpha}_2\\ v^{\alpha}_2 \end{smallmatrix}\right]\Big\|_{E^{\theta_2}_t(\alpha)}\Big)\\
&\times \Big(1+\Big\|\left[\begin{smallmatrix} u^{\alpha}_1\\ v^{\alpha}_1 \end{smallmatrix}\right]\Big\|^{\rho-1}_{E^{\theta_2}_t(\alpha)}
+\Big\|\left[\begin{smallmatrix} u^{\alpha}_2\\ v^{\alpha}_2 \end{smallmatrix}\right]\Big\|^{\rho-1}_{E^{\theta_2}_t(\alpha)}\Big),\nonumber
\end{split}
\]
for any $\left[\begin{smallmatrix} u^{\alpha}_1\\ v^{\alpha}_1 \end{smallmatrix}\right], \left[\begin{smallmatrix} u^{\alpha}_2\\ v^{\alpha}_2 \end{smallmatrix}\right] \in E^{\theta_2}_t(\alpha)$.\nonumber 
\end{lemma}

\proof The proof of \eqref{140r} was made in Bezerra et al. \cite[Lemma 3.2]{BCCN}. Given $\left[\begin{smallmatrix} u^{\alpha}_1\\ v^{\alpha}_1 \end{smallmatrix}\right], \left[\begin{smallmatrix} u^{\alpha}_2\\ v^{\alpha}_2 \end{smallmatrix}\right] \in E^{\theta_2}_t(\alpha)$, it follows from \cite[Lemma 3.1 and 3.2]{BCCN} and \eqref{hol-omega} that 
\[
\begin{split}
&\Big\|F\Big(t,\left[\begin{smallmatrix} u^{\alpha}_1\\ v^{\alpha}_1 \end{smallmatrix}\right]\Big)-F\Big(\tau,\left[\begin{smallmatrix} u^{\alpha}_2\\ v^{\alpha}_2 \end{smallmatrix}\right]\Big)\Big\|_{E^{\theta_1}_t(\alpha)}\\
&\leqslant\|f(u^{\alpha}_1)-f(u^{\alpha}_2)\|_{X^{\frac{\alpha\theta_1}{2}}}+\|\omega(t)(v^{\alpha}_1-v^{\alpha}_2)\|_{X^{\frac{\alpha\theta_1}{2}}}+\|(\omega(t)-\omega(\tau))v^{\alpha}_2\|_{X^{\frac{\alpha\theta_1}{2}}}\\
&\leqslant c\Big(\Big\|\left[\begin{smallmatrix} u^{\alpha}_1\\ v^{\alpha}_1 \end{smallmatrix}\right]-\left[\begin{smallmatrix} u^{\alpha}_2\\ v^{\alpha}_2 \end{smallmatrix}\right]\Big\|_{E^{\theta_2}_t(\alpha)}\Big)\Big(1+\Big\|\left[\begin{smallmatrix} u^{\alpha}_1\\ v^{\alpha}_1 \end{smallmatrix}\right]\Big\|^{\rho-1}_{E^{\theta_2}_t(\alpha)}+\Big\|\left[\begin{smallmatrix} u^{\alpha}_2\\ v^{\alpha}_2 \end{smallmatrix}\right]\Big\|^{\rho-1}_{E^{\theta_2}_t(\alpha)}\Big)+\kappa_0|t-\tau|^{\zeta}\Big\|\left[\begin{smallmatrix} u^{\alpha}_2\\ v^{\alpha}_2 \end{smallmatrix}\right]\Big\|_{E^{\theta_2}_t(\alpha)}\\
&\leqslant c(W,\kappa_0)\Big(|t-\tau|^{\zeta}+\Big\|\left[\begin{smallmatrix} u^{\alpha}_1\\ v^{\alpha}_1 \end{smallmatrix}\right]-\left[\begin{smallmatrix} u^{\alpha}_2\\ v^{\alpha}_2 \end{smallmatrix}\right]\Big\|_{E^{\theta_2}_t(\alpha)}\Big)\Big(1+\Big\|\left[\begin{smallmatrix} u^{\alpha}_1\\ v^{\alpha}_1 \end{smallmatrix}\right]\Big\|^{\rho-1}_{E^{\theta_2}_t(\alpha)}+\Big\|\left[\begin{smallmatrix} u^{\alpha}_2\\ v^{\alpha}_2 \end{smallmatrix}\right]\Big\|^{\rho-1}_{E^{\theta_2}_t(\alpha)}\Big),
\end{split}
\]
and the result is proved.
\qed

\medskip 

 In the next result, using \eqref{CondDissip} and exploiting gradient structure of \eqref{ffhghgf} we establish the global well posedness of \eqref{ffhghgf}.

\begin{theorem}\label{th:perturbed-problem}
 Suppose that \eqref{growth-condition}, \eqref{eq:supercritical} and \eqref{CondDissip} hold, and let $\omega$ be a decreasing strictly positive diffe\-rentiable function, $(\zeta,\kappa_0)-$H\"{o}lder continuous in $\mathbb{R}$ with \eqref{132r}, $\mu$ is a function satisfying \eqref{123} and \eqref{1506r}, and fix any number $s$ satisfying $\frac{N}2(1-\frac1\rho) -\frac1\rho<s<1$. For all $\alpha < 1$ close enough to $1$ the following statements hold:
\begin{itemize}
\item[$i)$] For any $\left[\begin{smallmatrix}u_{\tau\alpha}\\ v_{\tau\alpha} \end{smallmatrix}\right]\in X^\frac{s}2\times X^\frac{s-1}2$ the solution  $\left[\begin{smallmatrix}u^{\alpha}\\ v^{\alpha} \end{smallmatrix}\right]$ of \eqref{ffhghgf} obtained in Theorem \ref{th:local-well-posedness-of-perturbed-problem} exists globally in time and satisfies for $\xi>\tau\in\mathbb{R}$
\[
\Big\|\left[\begin{smallmatrix}u^{\alpha}(t)\\ v^{\alpha}(t) \end{smallmatrix}\right]\Big\|_{X^{\frac{1+\alpha}{4}}\times X^{\frac{-1+\alpha}{4}}}\leqslant C\big(\xi,\left[\begin{smallmatrix}u_{\tau\alpha}\\ v_{\tau\alpha} \end{smallmatrix}\right]\big),\quad \forall t\geqslant\xi,
\]
where $C$ is a positive constant which can be chosen uniformly for $\left[\begin{smallmatrix}u_{\tau\alpha}\\ v_{\tau\alpha} \end{smallmatrix}\right]$ in bounded subsets of $X^\frac{s}2\times X^\frac{s-1}2$.
\item[$ii)$] The family of maps
\begin{equation}\label{2126ro}
S_{\alpha}(t,\tau)\left[\begin{smallmatrix}u_{\tau\alpha}\\ v_{\tau\alpha} \end{smallmatrix}\right]=\left[\begin{smallmatrix} u^{\alpha}(t)\\ v^{\alpha}(t)\end{smallmatrix}\right],\quad\left[\begin{smallmatrix}u_{\tau\alpha}\\ v_{\tau\alpha} \end{smallmatrix}\right]\in X^\frac{s}2\times X^\frac{s-1}2,\quad t\geqslant\tau\in\mathbb{R}, 
\end{equation}
where $\left[\begin{smallmatrix}u^{\alpha}\\ v^{\alpha} \end{smallmatrix}\right]$ is a solution of \eqref{ffhghgf}. The nonlinear evolution process $\{S_{\alpha}(t,\tau):t\geqslant\tau\in\mathbb{R}\}$, is compact in $X^\frac{s}2\times X^\frac{s-1}2$.
\end{itemize}
\end{theorem}

\proof  
If $B$ is bounded in $X^\frac{s}2\times X^\frac{s-1}2$, using part $i)$ from Theorem \ref{th:local-well-posedness-of-perturbed-problem}, there exists a time $\tau_{B}>\tau\in\mathbb{R}$ such that for $\left[\begin{smallmatrix}u_{\tau\alpha}\\ v_{\tau\alpha} \end{smallmatrix}\right]\in B$, the solution $\left[\begin{smallmatrix} u^{\alpha}\\ v^{\alpha}\end{smallmatrix}\right]$ through $\left[\begin{smallmatrix}u_{\tau\alpha}\\ v_{\tau\alpha} \end{smallmatrix}\right]$ exists until $\tau_B$ and \eqref{eq:uniform-regularization} hold. Recalling Remark 3.3 $ i)$-$ii)$ and using \eqref{eq:uniform-regularization} and \eqref{equivaLR} we get that for each $\alpha< 1$ close enough to $1$ a constant $c_\alpha > 0$ exists such that
\begin{equation*}%\label{2156r}
\mathscr{L}_{\alpha}\Big(\left[\begin{smallmatrix}u^{\alpha}(\xi)\\ v^{\alpha}(\xi) \end{smallmatrix}\right]\Big)\leqslant c_{\alpha} \Big\|\left[\begin{smallmatrix}u^{\alpha}(\xi)\\ v^{\alpha}(\xi)\end{smallmatrix}\right]\Big\|_{X^{\frac{1+\alpha}{2}}\times X^{\frac{\alpha}{2}}}+D_3\leqslant c_{\alpha}M(\xi,B)+D_3,
\end{equation*}
since $X^\frac{1+\alpha}2\times X^\frac{\alpha}2\hookrightarrow X^{\frac{1+\alpha}{4}}\times X^{\frac{-1+\alpha}{4}}$.  On the other hand, we have from \eqref{12r} and \eqref{2216r} that 
\begin{equation*}
\mathscr{L}_{\alpha}\Big(\left[\begin{smallmatrix}u^{\alpha}(t)\\ v^{\alpha}(t) \end{smallmatrix}\right]\Big)\leqslant \mathscr{L}_{\alpha}\Big(\left[\begin{smallmatrix}u^{\alpha}(\xi)\\ v^{\alpha}(\xi) \end{smallmatrix}\right]\Big)+8D_{3}, \quad \xi\leqslant t<\xi_{u_{\tau\alpha}, v_{\tau\alpha}},
\end{equation*}
whereas from \eqref{equivaLR} we obtain
\begin{equation*}
D_1\Big\|\left[\begin{smallmatrix}u^{\alpha}(t)\\ v^{\alpha}(t)\end{smallmatrix}\right]\Big\|_{X^{\frac{1+\alpha}{4}}\times X^{\frac{-1+\alpha}{4}}}\leqslant\mathscr{L}_{\alpha}\Big(\left[\begin{smallmatrix}u^{\alpha}(\xi)\\ v^{\alpha}(\xi) \end{smallmatrix}\right]\Big)+9D_3\leqslant c_{\alpha}M(\xi,B)+10D_3,\quad \xi\leqslant t<\xi_{u_{\tau\alpha}, v_{\tau\alpha}}.
\end{equation*}
Since 
\begin{equation}\label{2246r}
X^{\frac{1+\alpha}{4}}\times X^{\frac{-1+\alpha}{4}}\hookrightarrow X^\frac{s}2\times X^\frac{s-1}2\text{ for all }\alpha\in(2s-1,1), 
\end{equation}
we now conclude all results of part $i)$.

Finally, part $ii)$ follows from part $i)$, because $S_{\alpha}(t,\tau)\left[\begin{smallmatrix}u_{\tau\alpha}\\ v_{\tau\alpha} \end{smallmatrix}\right]$ is bounded in $X^{\frac{1+\alpha}{4}}\times X^{\frac{-1+\alpha}{4}}$
and from compactness of the embedding \eqref{2246r}.
\qed

\medskip

For the global solutions obtained in Theorem \ref{th:perturbed-problem} we prove the existence of a pullback absorbing family.

\begin{theorem}\label{existence-absorvent}
There exists $R_\mathbb{A_\alpha}=C^{-1}_{1}(1+2M_1)$ such that the family $\mathbb{A}_{\alpha}= \{\mathbb{A}_{\alpha}(t): t\in\mathbb{R}\}$ is a pullback absorbing family in $X^{\frac{1+\alpha}{4}}\times X^{\frac{-1+\alpha}{4}}$ for the nonlinear evolution process $S_{\alpha}(\cdot,\cdot)$, where 
\[
\mathbb{A}_{\alpha}(t)=\Big\{\left[\begin{smallmatrix}u^{\alpha}\\ v^{\alpha}\end{smallmatrix}\right]\in X^{\frac{1+\alpha}{4}}\times X^{\frac{-1+\alpha}{4}}: \|\left[\begin{smallmatrix}u^{\alpha}\\ v^{\alpha}\end{smallmatrix}\right]\|_{X^{\frac{1+\alpha}{4}}\times X^{\frac{-1+\alpha}{4}}}\leqslant R_{\mathbb{A}_{\alpha}} \Big\}.
\]
\end{theorem}

\proof Let $\mathcal{B}_{\alpha}$ be a pullback-bounded family in $X^{\frac{s}{2}}\times X^{\frac{s-1}{2}}$ and, for $t\in\mathbb{R}$, let $$R(t)=\sup_{\tau\in(-\infty,t]}\mathscr{E}_{\mathcal{H}^\alpha_\tau}(\mathcal{B}_{\alpha}(\tau)),$$
which is finite for every $t$, due to the equivalence between the energy $\mathscr{E}_{\mathcal{H}^\alpha_t}$ and the $\mathcal{H}^\alpha_t-$norm. Using Lemma \ref{inequality-geral} and Lemma \ref{equiv-geral} for $(u_{\tau\alpha},v_{\tau\alpha})\in\mathscr{B}_{\alpha}(\tau)$ yields 
\begin{equation*}%\label{18r}
C_1\|\left[\begin{smallmatrix}u^{\alpha}\\ v^{\alpha}\end{smallmatrix}\right]\|_{X^{\frac{1+\alpha}{4}}\times X^{\frac{-1+\alpha}{4}}}\leqslant\mathscr{E}_{\mathcal{H}^\alpha_t}(u^\alpha,v^\alpha)\leqslant M_0R(t)e^{-\varepsilon_{\omega}(t)(t-\tau)}+M_1\leqslant 1+2M_1,
\end{equation*}
provided that $\tau\leqslant t_0=t_0(t):=t-\theta$, where
\[
\theta:=\max\Big\{0,\varepsilon^{-1}_{\omega}(t)\log\dfrac{M_0 R(t)}{1+M_1}\Big\}.
\]
Taking the supremum over $(u_{\tau\alpha},v_{\tau\alpha})\in\mathscr{B}_{\alpha}(\tau)$ and using \eqref{2126ro}, we obtain
\[
\|S_{\alpha}(t,\tau)\left[\begin{smallmatrix}u_{\tau\alpha}\\ v_{\tau\alpha} \end{smallmatrix}\right]\|_{X^{\frac{1+\alpha}{4}}\times X^{\frac{-1+\alpha}{4}}}=\|\left[\begin{smallmatrix}u^{\alpha}\\ v^{\alpha}\end{smallmatrix}\right]\|_{X^{\frac{1+\alpha}{4}}\times X^{\frac{-1+\alpha}{4}}}\leqslant R_{\mathbb{A}_{\alpha}}, \quad \forall \tau\leqslant t_0,
\]
which, reads exactly $S_{\alpha}(t,\tau)\mathscr{B}_{\alpha}(\tau)\subset \mathbb{A}_{\alpha}(t)$ whenever $\tau\leqslant t_0(t)$. Furthermore, $\mathbb{A}_{\alpha}$ is a pullback absorbing  in $X^{\frac{1+\alpha}{4}}\times X^{\frac{-1+\alpha}{4}}$ for the nonlinear evolution process $S_{\alpha}(\cdot,\cdot)$.
\qed

\medskip

Thanks to previous results we can prove a result of existence of time-dependent global attractor, with similar arguments used in Di Plinio, Duane, Temam \cite{DDT1,DDT2}. 

\begin{theorem}
Under the same conditions as in Theorem \ref{th:perturbed-problem}, there exists a unique pullback attractor $\mathscr{A}_{\alpha}=\{\mathscr{A}_{\alpha}(t):t\in\mathbb{R}\}$  in $X^{\frac{s}{2}}\times X^{\frac{s-1}{2}}$ for the nonlinear evolution process $\{S_{\alpha}(t,\tau):t\geqslant\tau\in\mathbb{R}\}$, where $\mathscr{A}_{\alpha}(t)=\omega_{\mathbb{A}_{\alpha}}(t)$.
\end{theorem}

\proof 
Let $\tau\in\mathbb{R}$ and $(u_{\tau\alpha},v_{\tau\alpha})\in\mathbb{A}_{\alpha}(\tau)$. Using \eqref{160r} and \eqref{Equivar}, we get
\begin{equation*}%\label{18r}
C_1\|\left[\begin{smallmatrix}u^{\alpha}\\ v^{\alpha}\end{smallmatrix}\right]\|_{X^{\frac{1+\alpha}{4}}\times X^{\frac{-1+\alpha}{4}}}\leqslant\mathscr{E}_{\mathcal{H}^\alpha_t}(u^\alpha,v^\alpha)\leqslant M_0C_2\|\left[\begin{smallmatrix}u_{\tau\alpha}\\ v_{\tau\alpha}\end{smallmatrix}\right]\|_{X^{\frac{1+\alpha}{4}}\times X^{\frac{-1+\alpha}{4}}}+M_1\leqslant M_0C_2 R_{\mathbb{A}_{\alpha}}+M_1,
\end{equation*}
provided that $\tau\leqslant t_0$. Thus, taking the supremum over $(u_{\tau\alpha},v_{\tau\alpha})\in\mathbb{A}_{\alpha}(\tau)$ and using \eqref{2126ro}, we concluded that $S_{\alpha}(t,\tau)\mathbb{A}_{\alpha}(\tau)$ is a bounded set in $X^{\frac{1+\alpha}{4}}\times X^{\frac{-1+\alpha}{4}}$ whenever $\tau\leqslant t_0(t)$. From  compactness of the embedding \eqref{2246r}, it follows that $S_{\alpha}(t,\tau)\mathbb{A}_{\alpha}(\tau)$ is a compact set in $X^{\frac{s}{2}}\times X^{\frac{s-1}{2}}$. Thus, for each $t\in\mathbb{R}$, then $\alpha_{t}(S_{\alpha}(t,\tau)\mathbb{A}_{\alpha}(\tau))=0$, for all $\tau\leqslant t_0(t)$. With this, \eqref{measurealpha} is satisfied. Therefore, the result follows by Theorems \ref{existence-attractor} and \ref{existence-absorvent} and Corollary \ref{uniqueness-attractor}. 
\qed

\medskip

However, under the same conditions as in Theorem \ref{th:perturbed-problem} we also can prove  a result of existence of pullback attractor  with similar arguments from Carvalho, Langa, Robinson \cite{CLR}; that is, the uniqueness is in the sense of part $2)$ of Remark \ref{remark:2.1} .

\begin{theorem}
Under the same conditions as in Theorem \ref{th:perturbed-problem}, there exists a unique pullback attractor $\mathscr{A}_{\alpha}=\{\mathscr{A}_{\alpha}(t):t\in\mathbb{R}\}$ for the nonlinear evolution process $\{S_{\alpha}(t,\tau):t\geqslant\tau\in\mathbb{R}\}$.
\end{theorem}

\proof 
It follows from Theorem \ref{existence-absorvent}, that there exists a pullback absorbing family $\mathbb{A}_{\alpha}$ which is bounded in $X^{\frac{1+\alpha}{4}}\times X^{\frac{-1+\alpha}{4}}$ whenever $\tau\leqslant t_0(t)$. Hence, from  compactness of the embedding \eqref{2246r}, it follows that $\mathbb{A}_{\alpha}$ is a pullback absorbing family of compact sets in $X^{\frac{s}{2}}\times X^{\frac{s-1}{2}}$. Therefore, the result follows from Carvalho, Langa, Robinson \cite[Theorem 2.12]{CLR}. 
\qed

\end{document}